\documentclass[11pt,reqno]{amsart}
\usepackage{amsmath,amssymb,amsthm,tabularx,color}
\usepackage{hyperref}
\usepackage[margin=25mm]{geometry}

\theoremstyle{plain}
\newtheorem{theorem}{Theorem}[section]
\newtheorem{proposition}[theorem]{Proposition}
\newtheorem{lemma}[theorem]{Lemma}
\newtheorem{corollary}[theorem]{Corollary}

\theoremstyle{definition}

\theoremstyle{remark}
\newtheorem{remark}[theorem]{Remark}
\newtheorem{example}[theorem]{Example}

\title[Boundary Representations \& Poisson Boundaries]{Boundary Representations \\ and Non-commutative Poisson Boundaries}

\author{Yoshimichi UEDA}
\address{
Graduate School of Mathematics, Nagoya University, 
Furocho, Chikusaku, Nagoya, 464-8602, Japan
}
\email{ueda@math.nagoya-u.ac.jp}

\date{\today}

\begin{document}

\begin{abstract}
As a sequel to \cite{HU2018}, we make explicit how Arveson's boundary representations and Izumi's non-commutative Poisson boundaries are related for the fixed point system of a unital completely positive (UCP) map: the Douglas and Coburn exact sequences are governed by the UCP idempotent obtained by averaging the iterates of the map, the quotient being the Poisson boundary at the von Neumann algebra level and the C$^*$-envelope, whose spectrum is the Choquet boundary, at the C$^*$-algebra level. The picture is illustrated by Toeplitz-type examples.
Moreover, it is remarked that the proof in \cite{HU2018} works even for unital $2$-positive maps, whereas Arveson's unitary equivalence theorem cannot be extended by replacing the assumption of complete order isomorphism with that of $2$-order isomorphism in accordance with the failure of Arveson's extension theorem beyond complete positivity.
\end{abstract}

\maketitle

\section{Introduction}

Boundary representations form the non-commutative Choquet boundary of an operator system \cite{Arv1969}, and it is known that there are always sufficiently many boundary representations; see \cite{Arv2008, DK2015} that crucially used the concept of maximal dilations for unital completely positive (UCP) maps due to \cite{DritschelMcCullough:JOT05}.  
On the other hand, the concept of non-commutative Poisson boundaries was introduced in Izumi's work \cite{Izumi2002} on product-type actions of the compact quantum group $\mathrm{SU}_q(2)$. Therefore, those two concepts whose names contain a common word ``boundary'' are of rather different characters. 

In \cite{HU2018}, a short proof of Arveson's boundary theorem was given in terms of non-commutative Poisson boundaries. More precisely, Theorem 1 of \cite{HU2018} asserts that, if $\varphi$ is a UCP selfmap on a unital C$^*$-algebra $B$ and $J \subseteq C^*(B^\varphi)$ is a simple essential ideal of $B$, then $\varphi=\mathrm{id}$ if and only if the quotient map $B\to B/J$ fails to be completely isometric on $B^\varphi$. The proof passes to the second dual $M=B^{**}$, uses the Ces\`aro--ultrafilter idempotent $\psi$ from $M$ onto $M^{\varphi}$, the Choi--Effros product $x\circ y=\psi(xy)$, and the fact that $M^\varphi\subseteq\mathrm{mult}(\psi)$, the multiplicative domain of $\psi$. The second half of \cite{HU2018} contains two observations that motivate us to provide this paper. First, Section 3.2 of \cite{HU2018} points out that the $*$-homomorphism $\rho=\psi\!\upharpoonright_{C^*(M^\varphi)}$ is, in the case where $M=B(H^2)$ and $\varphi=\mathrm{Ad}\,s^*$ with the unilateral shift $s$ on the Hardy space $H^2$, nothing but the Toeplitz symbol map, and more generally, that for a general normal UCP selfmap on a von Neumann algebra $M$ the sequence
\[
0\to\ker\rho \to C^*(M^\varphi)\to H^\infty(M,\varphi)\to0
\] 
is semisplit. This sequence is an analog of the Douglas sequence \cite[Theorem 7.11, Corollary 7.13]{Douglas1972} for Toeplitz operators with $L^\infty$-symbols, and thus $\rho$ may be viewed as an abstract generalization of the symbol map. Second, prompted by Davidson's (informal) proof via maximal dilations, the following ``vague'' question was raised in Section 3.3 of \cite{HU2018}: \emph{Is there any relation between the concepts of maximal dilations and non-commutative Poisson boundaries?} 

The purpose of this paper is to make the following picture explicit; it makes the two observations above precise and, we hope, gives a partial answer to the above question. In both the Coburn sequence 
\[
0\to\mathbb K\to\mathcal T\to C(\mathbb T)\to0
\]
(\cite{Coburn67}) and the Douglas sequence for $L^\infty$-symbols, the quotient is a boundary: the Choquet boundary in the first case (the spectrum of the $C^*$-envelope $C(\mathbb T)$) and the Poisson boundary $L^\infty(\mathbb T)$ in the second. Behind both there is a UCP map $\varphi$ (here $\mathrm{Ad}\,s^*$) and the UCP idempotent obtained by averaging its iterates: its range is the fixed point system, the Choi--Effros product turns the range into the boundary, and the symbol map is its restriction to the generated $C^*$-algebra. The measure-theoretic Poisson boundary always exists, whereas the topological Choquet boundary can be computed inside the $C^*$-algebra only when the idempotent lives there, that is, when $\varphi$ is weakly mean ergodic. Much of this seems implicit in the classical theory of Toeplitz operators, in Arveson's boundary theorem and in \cite{HU2018}; we make it explicit and show that it explains the known Toeplitz-type examples uniformly. 

The detailed contents are as follows. 

In Section 2, we will investigate the abstract symbol map as follows. The symbol map $\rho$ is shown to be independent of the choice of free ultrafilter, or more precisely, $\varphi^n(x) \to \rho(x)$ as $n\to\infty$ holds for every $x \in C^*(M^\varphi)$ without any use of ultrafilter or dilation (Proposition \ref{prop:sot}, but it is more or less Prunaru's work \cite{Pru11} and, for the unilateral shift, goes back to Barr\'{\i}a--Halmos work \cite{BH1982}). The kernel $\ker\rho$ is identified with the Shilov ideal, and hence $C^*_{\mathrm{env}}(M^\varphi)=H^\infty(M,\varphi)$ (Theorem \ref{thm:env}). Then, all the boundary representations for the operator system $M^\varphi$ are precisely compositions of the symbol map $\rho$ and  irreducible representations of the non-commutative Poisson boundary $H^\infty(M,\varphi)$, and they all are completely norming (Theorem \ref{thm:choquet}). This, in particular, proves that $M^\varphi$ has sufficiently many boundary representations without any use of the notion of maximal dilations. This may be regarded as a complete answer to the above-mentioned question posed in \cite{HU2018} in the case of fixed point systems of normal UCP selfmaps on von Neumann algebras. 

In Section 3, we will treat the $C^*$-algebra case; namely, we start with a UCP selfmap $\varphi$ on a unital $C^*$-algebra $B$. Our hope is that the idea of the previous section still works in the $C^*$-algebra case like \cite{HU2018}. There are two ways to apply our machinery to the pair $(B,\varphi)$. One is to follow the strategy of \cite{HU2018}, that is, one once passes to the second dual $M := B^{**}$ and take the canonical normal extension of $\varphi$ (denoted by the same symbol $\varphi$) to it. The other, motivated by Toeplitz operators with $L^\infty$-symbols, is to work in a specific representation of $B$, where $\varphi$ becomes normal. Our main target there is an operator system $S$ sitting in $B^\varphi$. Here, it is important to relax $S = B^\varphi$, because it is not trivial to find a UCP map whose fixed point system coincides with a given system. The boundary representations for $S$ are shown to be Poisson harmonic and, when they extend normally to $M$, correspond to type I direct summands of the non-commutative Poisson boundary $H^\infty(M,\varphi)$ (Theorem \ref{thm:bdry}). When $M=B^{**}$, by \cite{DK2015} we see that the resulting non-commutative Poisson boundary has sufficiently many (for $S$) type I direct summands. We then investigate the relationship between two approaches that $M = B^{**}$ and $M$ is given by the closure of $B$ in a concrete representation.  

In Section 4, we first investigate when the converse of Theorem \ref{thm:bdry} holds. A trivial situation where the converse fails explains that $S = B^\varphi$ is required. Then, when $M = B^{**}$, we show that the restriction $\rho_B$ of the symbol map $\rho : C^*(M^\varphi) \to H^\infty(M,\varphi)$ enjoys $\rho_B(B) \subseteq B$ if and only if $\varphi$ is \emph{weakly mean ergodic}. 
In a concrete representation, however, the condition $\rho_B(B) \subseteq B$ is strictly weaker: for the $C^*$-algebra $B$ generated by Toeplitz operators with $L^\infty$-symbols, $\rho_B(B)\subseteq B$ holds automatically while $\varphi$ is not mean ergodic (Remark~\ref{rem:Linfty}).  
Under these equivalent conditions, the key idea used in Section 2 works in the $C^*$-algebra situation, and the converse holds, and a Coburn-type semisplit exact sequence 
\[
0 \to J_S \to B \to C^*_\mathrm{env}(S) \to 0,   
\]
is obtained, and moreover, the ``Choquet boundary $\partial S$'' of $S$ is exactly the ``spectrum'' of $C^*_\mathrm{env}(S)$ (Theorem \ref{thm:conv}, Corollary \ref{cor:coburn}). Furthermore, in the same setting, the boundary representations for $S$ correspond to the minimal central projections of type I direct summands of $H^\infty(M,\varphi)$ with $M = B^{**}$ (Proposition \ref{prop:poisson}). It is also pointed out that weak$^*$ density of $C^*_{\mathrm{env}}(S)$ can fail in the setting of $M = B^{**}$ (Remark \ref{rem:notdense}). We then examine various examples: Toeplitz operators on the circle and the sphere, the Drury--Arveson shift, the Cuntz--Toeplitz algebra and the non-commutative disc algebra, boundary actions of free groups, Toeplitz--Pimsner algebras, and so on.   

In Appendix~\ref{app:2pos}, we show that the proof of \cite{HU2018} works for unital $2$-positive maps, whereas Arveson's extension theorem, the other pillar of the unitary equivalence criterion, does not survive beyond complete positivity \cite{CDPR}; we show that the criterion itself fails for $2$-order isomorphisms.

\bigskip\noindent
{\bf Use of AI tools:} This work was carried out with extensive use of an AI research assistant, Claude (Anthropic), mainly Sonnet models (versions not recorded), accessed via claude.ai, September--November 2026. The assistant contributed to nearly all stages of the work. The final manuscript was then rewritten by the author, who verified all proofs and cited statements, and subsequently polished with AI assistance. The author is solely responsible for the content of this paper. 

\medskip\noindent
{\bf Acknowledgements:} The author thanks his former student, Kei Hasegawa, for many discussions on various topics during his student days; Yasuhiko Sato for the comment mentioned at the beginning of Appendix~\ref{app:2pos}; and Kenta Kojin for discussions on operators arising in complex analysis.

\section{The non-commutative Poisson boundary and the C$^*$-envelope}
\label{sec:env}

Let $\varphi$ be a normal UCP selfmap on a von Neumann algebra $M$. Then we have an ultralimit
\[
\psi_\omega(x) := \lim_{m\to\omega}\frac{1}{m+1}\sum_{k=0}^{m} \varphi^k(x), \qquad x \in M
\]
(in the BW topology, see \cite{Arv1969}; that is, in the pointwise $\sigma$-weak topology, which is compact on bounded sets) with a free ultrafilter $\omega$ over the natural numbers $\mathbb{N}$, where $\varphi^k$ denotes the $k$-fold iterate of $\varphi$. 

This $\psi_\omega$ is a UCP idempotent of $M$ onto $M^\varphi$ with
\begin{equation}\label{varphi-invariance}
\psi_\omega = \psi_\omega\circ\varphi = \varphi\circ\psi_\omega.
\end{equation}
See e.g.\ \cite{HU2018}. It is important to remark that $\psi_\omega$ need not be normal. Nevertheless, it is well known that the Choi--Effros product $x\circ y := \psi_\omega(xy)$ (see \cite{CE1977}) makes the operator system $M^\varphi$ into a $W^*$-algebra $H^\infty(M,\varphi) = (M^\varphi,\circ)$, so that the identity map from $M^\varphi$ onto $H^\infty(M,\varphi)$ is a complete order isomorphism. It is known, see \cite[Section 2.5]{Izumi2002} (also \cite[Theorem 3.3]{Izumi2004}), that the $W^*$-algebra structure of $H^\infty(M,\varphi)$ does not depend on the choice of $\omega$, though $\psi_\omega$ may.

\begin{lemma}[{\cite[Claim 2]{HU2018}}]
\label{lem:mult}
Let $\mathrm{mult}(\psi_\omega)$ denote the multiplicative domain of the UCP idempotent $\psi_\omega$. Then
\[
M^\varphi \subseteq\mathrm{mult}(\psi_\omega), \qquad \varphi(\mathrm{mult}(\psi_\omega)) \subseteq \mathrm{mult}(\psi_\omega).
\]
In particular, the restriction $\rho:=\psi_\omega\!\upharpoonright_{C^*(M^\varphi)}$ is a unital $*$-homomorphism of $C^*(M^\varphi)$ onto $H^\infty(M,\varphi)$ with $\rho\!\upharpoonright_{M^\varphi}=\mathrm{id}$, which does not depend on the choice of the free ultrafilter $\omega$.
\end{lemma}
\begin{proof}
The two inclusions are explicitly stated in \cite[Claim 2]{HU2018} with proof. Since the multiplicative domain is a C$^*$-algebra, the first inclusion gives $C^*(M^\varphi)\subseteq\mathrm{mult}(\psi_\omega)$. Hence, for $a,b\in C^*(M^\varphi)$,
\[
\rho(a)\circ\rho(b)=\psi_\omega(\psi_\omega(a)\psi_\omega(b))=\psi_\omega(\psi_\omega(ab))=\rho(ab),
\]
so that $\rho$ is a unital $*$-homomorphism into $H^\infty(M,\varphi)$; it is onto because $\rho\!\upharpoonright_{M^\varphi}=\mathrm{id}$. Moreover,
\[
\rho(x_1^{\varepsilon_1}\cdots x_\ell^{\varepsilon_\ell}) = \psi_\omega(x_1^{\varepsilon_1}\cdots x_\ell^{\varepsilon_\ell}) = x_1^{\varepsilon_1}\circ \cdots \circ x_\ell^{\varepsilon_\ell}
\]
for any $x_i \in M^\varphi$ and $\varepsilon_i \in \{ \cdot, * \}$. Since the Choi--Effros product $x\circ y$ does not depend on $\omega$ (see the paragraph before this lemma; it also follows directly from Lemma~\ref{lem:monotone} below and polarization), and since the words $x_1^{\varepsilon_1}\cdots x_\ell^{\varepsilon_\ell}$ span a dense subspace of $C^*(M^\varphi)$, we conclude that $\rho$ is independent of $\omega$.
\end{proof}

In what follows, we call $\rho : C^*(M^\varphi) \to H^\infty(M,\varphi)$ the \emph{abstract symbol map}.

\begin{lemma}[Monotonicity] \label{lem:monotone}
For any $x\in M^\varphi$, we have $\varphi^k(x^*x) \nearrow \rho(x^*x)$ in $M$ as $k\to\infty$. Hence $d_x:=\rho(x^*x)-x^*x \geq 0$, $d_x\in C^*(M^\varphi)$, and $\varphi^k(d_x) \to 0$ strongly in any normal representation of $M$ as $k\to\infty$.
\end{lemma}
\begin{proof}
Let $x \in M^\varphi$ be fixed. By the Schwarz inequality, $x^* x = \varphi(x)^*\varphi(x) \leq \varphi(x^*x)$, and hence $\varphi^k(x^*x)$ is non-decreasing. Since this sequence is bounded, it converges to its supremum in $M$ in the $\sigma$-strong topology as $k\to\infty$. Since the Ces\`{a}ro means converge to the same limit, the supremum must be $\psi_\omega(x^*x)=\rho(x^*x) \in M^\varphi$ (in particular, it does not depend on $\omega$). Thus, $\varphi^k(d_x)=\rho(x^*x)-\varphi^k(x^*x) \to 0$ strongly in any normal representation of $M$ as $k\to\infty$.
\end{proof}

\begin{theorem}[$\ker\rho$ is the Shilov ideal]
\label{thm:env}
The kernel $\ker\rho$ is exactly the Shilov ideal $J_{M^\varphi}$ for $M^\varphi$ in $C^*(M^\varphi)$. Consequently, $\rho$ induces a $*$-isomorphism $C^*_\mathrm{env}(M^\varphi)=C^*(M^\varphi)/J_{M^\varphi} \cong H^\infty(M,\varphi)$.
\end{theorem}
\begin{proof}
Recall that $J_{M^\varphi}$ is the largest ideal $I$ of $C^*(M^\varphi)$ such that the quotient map is completely isometric on $M^\varphi$ (see \cite{Arv1969,Hamana1979}; the formulation is due to the former and the existence of such an ideal is due to the latter). The quotient map $C^*(M^\varphi)\to C^*(M^\varphi)/\ker\rho\cong H^\infty(M,\varphi)$ restricted to $M^\varphi$ is the identity map $M^\varphi\to H^\infty(M,\varphi)$, which is completely isometric. Thus, $\ker\rho$ is such an ideal, and hence $\ker\rho\subseteq J_{M^\varphi}$. Consequently, the quotient map $q : C^*(M^\varphi) \to C^*(M^\varphi)/J_{M^\varphi}$ factors through $H^\infty(M,\varphi)$ as $q=\sigma\circ\rho$ with a $*$-homomorphism $\sigma : H^\infty(M,\varphi)\to C^*(M^\varphi)/J_{M^\varphi}$. Since $M^\varphi=H^\infty(M,\varphi)$ as vector spaces and $\rho\!\upharpoonright_{M^\varphi}=\mathrm{id}$, we have $\sigma=q\!\upharpoonright_{M^\varphi}$, which is isometric. Hence $\sigma$ is isometric, and in particular injective, so that $J_{M^\varphi}=\ker q=\ker\rho$.
\end{proof}

\begin{remark}[Toeplitz example due to Arveson; see the appendix of \cite{Izu12}]
\label{rem:toeplitz}
Let $M=B(H^2)$, $\varphi(x)=s^* xs$ (with the unilateral shift $s$ on the Hardy space $H^2$). Then
$M^\varphi=\{T_f; f\in L^\infty(\mathbb T)\}$, $H^\infty(M,\varphi)\cong L^\infty(\mathbb T)$ and $\rho$ is the so-called symbol map of Douglas \cite{Douglas1972}, \cite[\S3.2]{HU2018}; Theorem~\ref{thm:env} says that its kernel is the Shilov ideal of the Toeplitz operator system.
\end{remark}

\begin{theorem}[Boundary representations for $M^\varphi$ are irreducible representations of the non-commutative Poisson boundary]
\label{thm:choquet}
Let $\varphi$ be a normal UCP selfmap on a von Neumann algebra $M$. Here $H^\infty(M,\varphi)$ is regarded as a C$^*$-algebra, and its representations are not assumed to be normal (in the Toeplitz example of Remark~\ref{rem:toeplitz}, they include the characters of $L^\infty(\mathbb T)$, which are in general not normal). For every unital representation $\pi$ of $H^\infty(M,\varphi)$, the representation $\pi\circ\rho$ of $C^*(M^\varphi)$ is the \emph{unique} UCP extension of its restriction to $M^\varphi$. Consequently, we have:
\begin{itemize}
\item[(1)] An irreducible representation of $C^*(M^\varphi)$ is a boundary representation for $M^\varphi$ if and only if it annihilates $\ker\rho=J_{M^\varphi}$; the boundary representations for $M^\varphi$ are exactly $\pi\circ\rho$ with irreducible representations $\pi$ of $H^\infty(M,\varphi)$. In other words, $\pi\mapsto\pi\circ\rho$ induces a bijection from the spectrum of $H^\infty(M,\varphi)$ (the set of unitary equivalence classes of its irreducible representations) onto the Choquet boundary of $M^\varphi$ (the set of unitary equivalence classes of its boundary representations).
\item[(2)] The boundary representations for $M^\varphi$ are completely norming for $M^\varphi$.
\end{itemize}
\end{theorem}
\begin{proof}
Let $\theta : C^*(M^\varphi)\to B(\mathcal{K})$ be a UCP map with $\theta\!\upharpoonright_{M^\varphi} = \pi\circ\rho\!\upharpoonright_{M^\varphi}$. By Lemma~\ref{lem:monotone} we have, for every $x \in M^\varphi$,
\[
\theta(x^*x)\leq \theta(\rho(x^* x))=\theta(x^*\circ x)=\pi(x^*\circ x)=\pi(x)^*\pi(x)=\theta(x)^*\theta(x).
\]
Now the Schwarz inequality gives the reverse inequality; namely, $\theta(x^*x)=\theta(x)^*\theta(x)$. Applying this to $x^*\in M^\varphi$ as well, we conclude that $M^\varphi \subseteq \mathrm{mult}(\theta)$, implying that $C^*(M^\varphi) \subseteq \mathrm{mult}(\theta)$. Thus, $\theta$ is a $*$-homomorphism that agrees with $\pi\circ\rho$ on $M^\varphi$, and hence on $C^*(M^\varphi)$. Therefore, $\pi\circ\rho=\theta$ is the unique UCP extension of its restriction to $M^\varphi$.

(1) If $\pi$ is an irreducible representation of $H^\infty(M,\varphi)$, then $\pi\circ\rho$ is an irreducible representation of $C^*(M^\varphi)$ (since $\rho$ is surjective), and it is a boundary representation by the first part. Conversely, Theorem~\ref{thm:env} asserts that any boundary representation annihilates the Shilov ideal $J_{M^\varphi}=\ker\rho$, and hence it is of the form $\pi\circ\rho$ with an irreducible representation $\pi$ of $H^\infty(M,\varphi)$, since $\rho$ is surjective. The same argument shows that an irreducible representation annihilating $\ker\rho$ is a boundary representation. Finally, if $\pi\circ\rho$ and $\pi'\circ\rho$ are unitarily equivalent through a unitary $U$, then $U\pi(\rho(a))U^*=\pi'(\rho(a))$ for all $a\in C^*(M^\varphi)$, and hence $\pi$ and $\pi'$ are unitarily equivalent through the same $U$, since $\rho$ is surjective; this gives the injectivity of the induced map on equivalence classes.

(2) $H^\infty(M,\varphi)$ is a C$^*$-algebra, so its irreducible representations are completely norming for it, and $\rho\!\upharpoonright_{M^\varphi}$ is the identity map of $M^\varphi$ onto $H^\infty(M,\varphi)$, which is completely isometric. Hence, for $X\in M_n(M^\varphi)$, we have $\Vert X\Vert=\sup_\pi\Vert(\pi\circ\rho)_n(X)\Vert$, where $\pi$ runs over all irreducible representations of $H^\infty(M,\varphi)$.
\end{proof}

When $M$ is finite dimensional, $H^\infty(M,\varphi)=(M^\varphi,\circ)$ is a finite dimensional C$^*$-algebra, and Theorems~\ref{thm:env} and~\ref{thm:choquet} say that it is the C$^*$-envelope of $M^\varphi$ and that the boundary representations are exactly its irreducible representations.

We remark that item (2) above was established without using the work of Davidson--Kennedy \cite{DK2015}. Thus, in this special case, the non-commutative Poisson boundary alone yields sufficiently many boundary representations.

\medskip
For the reader's convenience, we will reconstruct Prunaru's work \cite{Pru11} in our terminology, which can be understood as a generalization of Barr\'{\i}a--Halmos's work \cite{BH1982} on the symbol map for Toeplitz operators in our context. 

\begin{proposition}[Prunaru {\cite[Theorem 2]{Pru11}}]
\label{prop:sot}
For every $x\in C^*(M^\varphi)$ the sequence $\varphi^k(x)$ converges strongly to $\rho(x) \in M^\varphi$ (in any normal representation of $M$) as $k\to\infty$. In particular, no ultrafilter, averaging or dilation is needed to construct the abstract symbol map $\rho$.
\end{proposition}
\begin{proof}
We use the same notation as in Lemma~\ref{lem:monotone}.

For any $x \in M^\varphi$ and $a \in C^*(M^\varphi)$ we observe that
\[
\rho(xa)-xa = x\circ\rho(a) - xa = (\rho(x\rho(a)) - x\rho(a)) + x(\rho(a)-a).
\]
Using this repeatedly (together with the density of the linear span of the words in elements of $M^\varphi$ in $C^*(M^\varphi)$) we obtain that $\ker\rho$ is the norm-closed left ideal generated by the $\rho(xy)-xy$ with $x,y \in M^\varphi$. By the polarization trick together with the fact that $M^\varphi$ is a self-adjoint subspace, we see that $\ker\rho$ is the norm-closed left ideal generated by $d_x:=\rho(x^* x) -x^*x$ with $x \in M^\varphi$.

We may and do assume $M \subset B(\mathcal{H})$ with a Hilbert space $\mathcal{H}$. For any $a \in C^*(M^\varphi)$, $x \in M^\varphi$ and $\xi \in \mathcal{H}$, we have, by the Schwarz inequality and Lemma~\ref{lem:monotone},
\[
\Vert\varphi^k(ad_x^{1/2})\xi\Vert^2 \leq (\varphi^k(d_x^{1/2}a^*a\,d_x^{1/2})\xi|\xi) \leq \Vert a\Vert^2(\varphi^k(d_x)\xi|\xi) \to 0
\]
as $k\to \infty$. Since $C^*(M^\varphi)d_x\subseteq C^*(M^\varphi)d_x^{1/2}$ and $\Vert\varphi^k\Vert \leq 1$ (as a linear map), we conclude that $\varphi^k(b) \to 0$ strongly for every $b\in \ker\rho$. The desired assertion immediately follows, since $x-\rho(x)\in\ker\rho$ and $\varphi^k(\rho(x))=\rho(x)$ for $x\in C^*(M^\varphi)$.
\end{proof}

\begin{remark}
Arveson's observation, written by Izumi \cite[Appendix]{Izu12}, gives a dilation picture of the above fact and even more in our framework. Assume that $M$ has separable predual, and take the minimal dilation $(K,p,\alpha)$ ($M=pKp$, $\varphi^k(a)=p\alpha^k(a)p$ for all $k \in \mathbb{N}$ with a $*$-automorphism $\alpha$ of $K$). Then, the map
$x\mapsto pxp$ is a bijective $*$-homomorphism from $K^\alpha$ onto $H^\infty(M,\varphi)$, and $\varphi^k(xy)\to x\circ y = \rho(xy)$ strongly for any $x,y \in M^\varphi$; see \cite[Corollary 5.2]{Izu12}. In the Toeplitz example, we have $K=B(L^2(\mathbb T))$, $\alpha=\mathrm{Ad}\,u^*$ with the bilateral shift $u$, $K^\alpha=L^\infty(\mathbb T)$ and $p$ is the projection onto $H^2$.
\end{remark}

\section{Boundary representations inside the non-commutative Poisson boundary}
\label{sec:bdry}

We start with C$^*$-algebraic data consisting of a unital C$^*$-algebra $B$, a UCP selfmap $\varphi$ on $B$ and an operator system $S \subseteq B^\varphi$ such that $B=C^*(S)$.

\medskip
Throughout this section, we work in the following situation: there are a von Neumann algebra $M$ and a normal UCP selfmap on $M$ extending $\varphi$ such that $B$ is a $\sigma$-weakly dense unital C$^*$-subalgebra of $M$. We denote this normal extension by the same letter $\varphi$; thus $\varphi(B) \subseteq B$. Note that $S\subseteq B^\varphi=B\cap M^\varphi$, and hence $B=C^*(S)\subseteq C^*(M^\varphi)$. We call such a pair $(M,\varphi)$ a \emph{$W^*$-envelope} for the pair $(B,\varphi)$ in what follows. 
Note that the assumption $B=C^*(S)$ is used only to ensure $B\subseteq C^*(M^\varphi)$, on which $\rho=\psi_\omega\!\upharpoonright_{C^*(M^\varphi)}$ is a $*$-homomorphism.

\medskip
A $W^*$-envelope $(M,\varphi)$ always exists: as in \cite{HU2018}, take the second dual $M:=B^{**}$ with the second adjoint of $\varphi$ (denoted by $\varphi^{**}$ when we want to emphasize it), which is a normal UCP map since $M_n(B)^{**}=M_n(B^{**})$; see e.g.\ the proof of \cite[Theorem 2.1.1]{Arv1972}. Other choices of $(M,\varphi)$ are provided by concrete representations of $B$ (see Proposition~\ref{prop:bridge}). 

\medskip
Then, we have the non-commutative Poisson boundary $H^\infty(M,\varphi)$ and the abstract symbol map $\rho : C^*(M^\varphi) \to H^\infty(M,\varphi)$, which is the restriction of a UCP idempotent $\psi_\omega : M \to M^\varphi$ for any free ultrafilter $\omega$ over $\mathbb{N}$. We denote its restriction to $B=C^*(S)\subseteq C^*(M^\varphi)$ by $\rho_B : B \to H^\infty(M,\varphi)$, which is a unital $*$-homomorphism.

\medskip
We now recall the following well-known facts concerning the second dual. Let $\pi : B \curvearrowright \mathcal{H}_\pi$ be a representation. Then, there is a unique normal representation $\bar{\pi}$ of the second dual $B^{**}$ on the same $\mathcal{H}_\pi$ extending $\pi$ along $B \subset B^{**}$, and $\bar{\pi}(B^{**}) = \pi(B)''$ holds. More generally, any unital $*$-homomorphism from $B$ into a von Neumann algebra extends uniquely to a normal $*$-homomorphism from $B^{**}$ into the same von Neumann algebra. See e.g.\ \cite[Section III.5.2.9]{Bla2006}. 
\emph{These facts concern $M=B^{**}$.} For a general $M$ as above, a representation of $B$ need not extend to a normal representation of $M$ (for instance, if $B$ is the Toeplitz algebra and $M=B(H^2)$, the characters of $B$ coming from the symbol map do not); when it does, the extension is unique, and we use the same notation $\bar{\pi}$.

When $M=B^{**}$, applying the above to $\rho_B : B \hookrightarrow C^*(M^\varphi) \to H^\infty(M,\varphi)$ one has a unique normal $*$-homomorphism $\bar{\rho}_B : M \to H^\infty(M,\varphi)$ extending $\rho_B$, and its range is the $\sigma$-weak closure of $\rho_B(B)$ in $H^\infty(M,\varphi)$. One cannot conclude that $\bar{\rho}_B(M) = H^\infty(M,\varphi)$ in general. 

\medskip
In what follows, a representation $\pi : B \curvearrowright \mathcal{H}_\pi$ is said to be \emph{$\varphi$-harmonic} if $\pi\circ\varphi = \pi$ holds on $B$.

\begin{proposition}[Harmonic factorization]
\label{prop:harm}
Let $\pi$ be a unital representation of $B$ with a (unique) normal extension $\bar{\pi}$ to $M$. The following are equivalent:
\begin{itemize}
\item[(a)] $\pi\circ\varphi=\pi$ on $B$, that is, $\pi$ is $\varphi$-harmonic.
\item[(b)] $\bar{\pi}\circ\varphi=\bar{\pi}$ on $M$.
\item[(c)] $\bar{\pi}\circ\psi_\omega=\bar{\pi}$ on $M$ for some free ultrafilter $\omega$.
\item[(d)] $\bar{\pi}\circ\psi_\omega=\bar{\pi}$ on $M$ for every free ultrafilter $\omega$.
\end{itemize}
In this case, $\hat\pi:=\bar{\pi}\!\upharpoonright_{M^\varphi}$ is a normal representation of $H^\infty(M,\varphi)$ and $\pi = \hat\pi\circ\rho_B$ holds.
\end{proposition}
\begin{proof}
(a)$\Leftrightarrow$(b): The maps $\bar{\pi}\circ\varphi$ and $\bar\pi$ are normal on $M$ and agree on $B$ by (a), and $B$ is $\sigma$-weakly dense in $M$ (by the standing assumption; if $M=B^{**}$, this is Goldstine's theorem).

(b)$\Rightarrow$(d): Observe that
\[
\bar{\pi}\left(\frac{1}{m+1}\sum_{k=0}^m\varphi^k(x)\right)=\bar{\pi}(x), \qquad x \in M,
\]
and take the $\omega$-limit for any free ultrafilter $\omega$, since $\bar{\pi}$ is normal.

(d)$\Rightarrow$(c) is trivial.

(c)$\Rightarrow$(b): By \eqref{varphi-invariance}, we have $\bar{\pi}\circ\varphi=\bar{\pi}\circ\psi_\omega\circ\varphi=\bar{\pi}\circ\psi_\omega=\bar{\pi}$ on $M$.

For any $x,y \in M^\varphi$ we also have $\bar{\pi}(x\circ y)=\bar{\pi}(\psi_\omega(xy))=\bar{\pi}(xy)=\bar{\pi}(x)\bar{\pi}(y)$ by (c), so that $\hat\pi$ is a unital $*$-representation of $H^\infty(M,\varphi)$. Since the identity map $M^\varphi\to H^\infty(M,\varphi)$ is a complete order isomorphism, the $\sigma$-weak topology of $H^\infty(M,\varphi)$, which is determined by its order structure, coincides with the relative $\sigma$-weak topology of $M^\varphi\subseteq M$; hence $\hat\pi$ is normal because $\bar\pi$ is. 
Finally, $\hat\pi\circ\rho_B=\bar\pi\circ\psi_\omega\!\upharpoonright_B=\bar\pi\!\upharpoonright_B=\pi$ by (c). 
\end{proof}

\begin{theorem}[Boundary representations are Poisson harmonic]
\label{thm:bdry}
Let $\pi$ be a boundary representation for $S$, that is, $\pi$ is an irreducible representation of $B=C^*(S)$ on a Hilbert space $\mathcal{H}_\pi$ such that $\pi$ is the only UCP extension of the restriction $\pi\!\upharpoonright_S$ to $B$. Then:
\begin{itemize}
\item[(1)] $\pi$ is $\varphi$-harmonic.
\item[(2)] If $\pi$ has a (unique) normal extension $\bar{\pi}$ to $M$, then there is a unique central projection $z_\pi$ of $H^\infty(M,\varphi)$ such that $z_\pi\circ H^\infty(M,\varphi) \cong B(\mathcal{H}_\pi)$, a type I factor direct summand. Inequivalent boundary representations give mutually orthogonal $z_\pi$.
\end{itemize}
When $M = B^{**}$, we also have
\begin{itemize}
\item[(3)] Let $z_\mathrm{bdy}$ be the join of all the $z_\pi$ with boundary representations $\pi$ for $S$ in the projection lattice of the center $\mathcal{Z}(H^\infty(M,\varphi))$. Then, the mapping $a \in B \mapsto z_\mathrm{bdy}\circ \rho_B(a) \in H^\infty(M,\varphi)$ is completely isometric on $S$.
\end{itemize}
\end{theorem}
\begin{proof}
(1) Since $\pi\circ\varphi$ is a UCP extension of $\pi\!\upharpoonright_S$ to $B$ (recall $\varphi\!\upharpoonright_S=\mathrm{id}$), we have $\pi\circ\varphi=\pi$ by the unique extension property.

\medskip
(2) By (1) and Proposition~\ref{prop:harm}, $\hat\pi=\bar\pi\!\upharpoonright_{M^\varphi}$ is a normal representation of $H^\infty(M,\varphi)$ with $\pi=\hat\pi\circ\rho_B$. Since $\pi(B)'' = B(\mathcal{H}_\pi)$ (by the irreducibility of $\pi$) and $\pi(B) = \hat\pi(\rho_B(B)) \subseteq \hat\pi(H^\infty(M,\varphi))$, and since $\hat\pi(H^\infty(M,\varphi))$ is a $\sigma$-weakly closed $*$-subalgebra of $B(\mathcal H_\pi)$ (as $\hat\pi$ is normal), we conclude that $\hat{\pi} : H^\infty(M,\varphi) \to B(\mathcal{H}_\pi)$ is surjective. By a standard fact on von Neumann algebras (a normal surjective $*$-homomorphism has a central support), there is a unique central projection $z_\pi$ of $H^\infty(M,\varphi)$ such that $\hat{\pi}$ is injective on $z_\pi\circ H^\infty(M,\varphi)$ and annihilates $(1-z_\pi)\circ H^\infty(M,\varphi)$. Hence $z_\pi \circ H^\infty(M,\varphi) \cong B(\mathcal{H}_\pi)$ by $z_\pi \circ x \mapsto \hat{\pi}(x)$.

Let $\pi_1,\pi_2$ be inequivalent boundary representations for $S$ having normal extensions. Then $\hat\pi_1,\hat\pi_2$ are irreducible normal representations of $H^\infty(M,\varphi)$, and they must be inequivalent: otherwise $\pi_i=\hat\pi_i\circ\rho_B$, $i=1,2$, would be equivalent. Therefore, their central supports $z_{\pi_1},z_{\pi_2}$ must be orthogonal to each other.

\medskip
(3) Let $M=B^{**}$, so that every boundary representation has a normal extension. By \cite{Arv2008} and \cite{DK2015}, there are sufficiently many boundary representations for $S$; namely, the direct sum of all of them is completely isometric on $S$. For $X\in M_n(S)$ and a boundary representation $\pi$ we have $\pi_n(X)=\hat\pi_n(\rho_B(X))=\hat\pi_n(z_\pi\circ\rho_B(X))$, so that $\Vert\pi_n(X)\Vert\leq\Vert(z_\mathrm{bdy}\circ\rho_B)_n(X)\Vert\leq\Vert X\Vert$ (the last inequality because $\rho_B$ is a $*$-homomorphism and $z_{\rm bdy}\circ(\cdot)$ is a $*$-homomorphism onto a corner). Taking the supremum over all boundary representations $\pi$ gives $\Vert X\Vert=\Vert(z_\mathrm{bdy}\circ\rho_B)_n(X)\Vert$. 
\end{proof}

In contrast to Theorem~\ref{thm:choquet}(2), item (3) above heavily depends upon \cite{Arv2008,DK2015}. Giving a proof of it without using such a heavy result is an interesting effort from the viewpoint of the question posed in \cite{HU2018}.

\medskip
In closing this section, we prove a proposition showing two special properties of the non-commutative Poisson boundary for the second dual $W^*$-envelope $(B^{**},\varphi^{**})$.

\begin{proposition}
\label{prop:bridge}
The following hold for $H^\infty(B^{**},\varphi^{**})$:
\begin{itemize}
\item[(1)] There is a unique central projection $z$ of $B^{**}$ such that $\bar{\rho}_B$ restricts to an isomorphism of $zB^{**}$ onto $\bar{\rho}_B(B^{**}) \subseteq H^\infty(B^{**},\varphi^{**})$.
\item[(2)] Let $N$ be a von Neumann algebra with a normal UCP selfmap $\theta$ on $N$, and let $\sigma : B \to N$ be a unital $*$-homomorphism such that $\sigma(B)'' = N$ and $\theta\circ\sigma = \sigma\circ\varphi$. Then, the normal extension $\bar{\sigma} : B^{**} \to N$ induces a normal $*$-homomorphism from $H^\infty(B^{**},\varphi^{**})$ onto $H^\infty(N,\theta)$; thus $H^\infty(N,\theta)$ can be realized as a central summand of $H^\infty(B^{**},\varphi^{**})$.
\end{itemize}
\end{proposition}
\begin{proof} 
(1) The kernel of the normal $*$-homomorphism $\bar\rho_B$ is a $\sigma$-weakly closed ideal of the von Neumann algebra $B^{**}$, and hence of the form $(1-z)B^{**}$ with a unique central projection $z$.

(2) The maps $\bar\sigma\circ\varphi^{**}$ and $\theta\circ\bar\sigma$ are normal and agree on $B$, because $\bar\sigma\circ\varphi^{**}\!\upharpoonright_B=\sigma\circ\varphi=\theta\circ\sigma=\theta\circ\bar\sigma\!\upharpoonright_B$. Therefore, they coincide on $B^{**}$. Since $\bar\sigma$ is normal, it commutes with the Ces\`aro--ultrafilter limits: $\bar\sigma\circ\psi_\omega=\psi_{\omega,N}\circ\bar\sigma$, where $\psi_{\omega,N}$ denotes the corresponding UCP idempotent of $N$ onto $N^\theta$. In particular, $\bar\sigma((B^{**})^{\varphi^{**}})\subseteq N^\theta$. Moreover, $\bar\sigma(B^{**})=\sigma(B)''=N$. Thus, if $y\in N^\theta$ and $y=\bar\sigma(x)$, then $y=\psi_{\omega,N}(y)=\bar\sigma(\psi_\omega(x))$ with $\psi_\omega(x)\in(B^{**})^{\varphi^{**}}$; thus $\bar\sigma$ maps $(B^{**})^{\varphi^{**}}$ onto $N^\theta$. Finally, for any $x,y\in(B^{**})^{\varphi^{**}}$,
\[
\bar\sigma(x\circ y)=\bar\sigma(\psi_\omega(xy))=\psi_{\omega,N}(\bar\sigma(x)\bar\sigma(y))=\bar\sigma(x)\circ\bar\sigma(y).
\]
Thus the restriction of $\bar\sigma$ is a $*$-homomorphism of $H^\infty(B^{**},\varphi^{**})$ onto $H^\infty(N,\theta)$, which is normal since the $\sigma$-weak topologies of both non-commutative Poisson boundaries are the relative ones (see the proof of Proposition~\ref{prop:harm}). Its kernel is a $\sigma$-weakly closed ideal, a corner $(1-z')\circ H^\infty(B^{**},\varphi^{**})$ with a central projection $z'$, and $H^\infty(N,\theta)\cong z'\circ H^\infty(B^{**},\varphi^{**})$.
\end{proof}

\section{$\varphi$-harmonicity implies being a boundary representation}
\label{sec:ex}
 
In contrast to Section~\ref{sec:env}, the converse of Theorem~\ref{thm:bdry}, i.e., the implication ($\varphi$-harmonicity $\Rightarrow$ being a boundary representation), is not automatic in the $C^*$-algebra setting. We investigate when the converse holds, and describe what follows from it (an exact sequence and the position of the non-commutative Poisson boundary).
 
\medskip
We start with two examples.
 
\begin{example}[Toeplitz and the choice of $M$]
In Remark~\ref{rem:toeplitz} the algebra is $N=B(H^2)$, not a second dual. Let $B=C^*(s)$ be the Toeplitz algebra and $\sigma\colon B\hookrightarrow B(H^2)$ the inclusion, so that $\sigma(B)''=N$, and let $\theta=\mathrm{Ad}\,s^*$ on $N$ (then $\theta\circ\sigma=\sigma\circ\varphi$ with $\varphi=\mathrm{Ad}\,s^*\!\upharpoonright_B$).
By Proposition~\ref{prop:bridge}(2), $L^\infty(\mathbb T)=H^\infty(N,\theta)$ is a corner of the non-commutative Poisson boundary $H^\infty(B^{**},\varphi^{**})$ computed in the second dual; we do not claim the two coincide.
\end{example}
 
The example tells us that there are several candidates of $W^*$-envelopes for $(B,\varphi)$, and the second dual $M=B^{**}$ gives a universal framework.
 
\begin{example}[The converse of Theorem~\ref{thm:bdry} fails in general]
Take $B=C(\bar{\mathbb D})$, $\varphi=\mathrm{id}$ and $S=\mathrm{span}\{1,z,\bar z\}$, which is an operator system with $C^*(S)=B$ (Stone--Weierstrass) and $\varphi\!\upharpoonright_S=\mathrm{id}$.
Every point evaluation $\delta_x$ is harmonic.
The unique representing measure of $\delta_x\!\upharpoonright_S$ is $\delta_x$ if and only if $x$ is an extreme point of $\bar{\mathbb D}$, i.e., $|x|=1$.
Thus, $\delta_x$ with $|x|<1$ is a harmonic irreducible representation which is not a boundary representation of $S$.
Here $S\subsetneqq B^\varphi=B$, and hence the meaningful question is the case $S=B^\varphi$ (see the discussion after Theorem~\ref{thm:conv}).
\end{example}
 
By this example we concentrate on the case $S=B^\varphi$ in what follows.
 
\medskip
We have constructed the $*$-homomorphism $\rho_B:B\to H^\infty(M,\varphi)$, but it is non-trivial whether or not $\rho_B(B)=\psi_\omega(B)\subseteq B^\varphi\subseteq M^\varphi$ holds. If this is the case, it is natural to believe that the machinery in Section~\ref{sec:env} works even in the $C^*$-algebra case.
 
\medskip
Here is a technical lemma, which is probably well known (see e.g., \cite[Section 7.2]{EF2025} and historical comments therein).
 
\begin{lemma}[Mean ergodicity gives $\psi$]
\label{lem:meanergodic}
Let $\varphi$ be a UCP selfmap on a unital C$^*$-algebra $B$, and set
\[
\mathcal{C}_m(b):=\frac{1}{m+1}\sum_{k=0}^{m}\varphi^k(b), \qquad b\in B.
\]
Assume that, for every $b\in B$, the sequence $\mathcal{C}_m(b)$ has a weak cluster point in $B$ (equivalently, some subnet of it converges weakly). Then $\mathcal{C}_m(b)$ converges in norm for every $b\in B$, and $\psi(b):=\lim_{m\to\infty}\mathcal{C}_m(b)$, $b\in B$, is a UCP idempotent of $B$ onto $B^\varphi$ with $\psi\circ\varphi=\varphi\circ\psi=\psi$.
Moreover, $\pi\circ\psi=\pi$ for every representation $\pi$ of $B$ with $\pi\circ\varphi=\pi$.
\end{lemma}
\begin{proof}
Let $y$ be a weak cluster point of the sequence $\mathcal{C}_m(b)$. Since $\varphi$ is weakly continuous and $\varphi(\mathcal{C}_m(b))-\mathcal{C}_m(b)=\frac{1}{m+1}(\varphi^{m+1}(b)-b)\to0$ in norm, we get $\varphi(y)=y$.
Since
\[
b-\mathcal{C}_m(b)=\frac{1}{m+1}\sum_{k=0}^{m}(b-\varphi^k(b))\in\mathrm{ran}(\mathrm{id}-\varphi),
\]
$b-y$ lies in the weak closure of the subspace $\mathrm{ran}(\mathrm{id}-\varphi)$, which equals its norm closure.
For any $z=(\mathrm{id}-\varphi)u$ with $u\in B$ we have $\mathcal{C}_m(z)=\frac{1}{m+1}(u-\varphi^{m+1}(u))\to0$ in norm; since $\Vert\mathcal{C}_m\Vert\leq1$, $\mathcal{C}_m(w)\to0$ in norm for every $w$ in the norm closure of $\mathrm{ran}(\mathrm{id}-\varphi)$. Hence $\mathcal{C}_m(b-y)\to0$, and $\mathcal{C}_m(y)=y$ gives $\mathcal{C}_m(b)\to y$ in norm. In particular $y$ is unique, and we define $\psi(b):=y$.
Note that $\psi$ is a norm limit of the UCP maps $\mathcal{C}_m$, and hence $\psi$ is UCP. Moreover, $\psi\circ\varphi=\psi$ follows from $\mathcal{C}_m(\varphi(b))-\mathcal{C}_m(b)\to0$ and $\varphi\circ\psi=\psi$ from $\varphi(y)=y$. For every $b\in B^\varphi$, we have $\mathcal{C}_m(b)=b$, and thus $\psi$ is an idempotent onto $B^\varphi$. Finally, if $\pi\circ\varphi=\pi$, then $\pi\circ\mathcal{C}_m=\pi$, and hence $\pi\circ\psi=\pi$ by the norm continuity of $\pi$.
\end{proof}
 
We say that a UCP selfmap $\varphi$ on a unital $C^*$-algebra $B$ is \emph{weakly mean ergodic} if $\varphi$ fulfills the assumption of Lemma~\ref{lem:meanergodic}. By the lemma, this is equivalent to the norm convergence of $\mathcal{C}_m(b)$ for every $b\in B$, and in this case $\psi_\omega\!\upharpoonright_B=\psi$ for every free ultrafilter $\omega$ (norm convergence implies weak$^*$ convergence in $M$).
 
\medskip
We show that $\rho_B(B)\subseteq B$ is actually equivalent to the weak mean ergodicity of $\varphi$.
 
\begin{proposition}
\label{prop:rho_B}
When $M=B^{**}$, the following are equivalent:
\begin{itemize}
\item[(a)] $\rho_B(B)\subseteq B$.
\item[(b)] $\varphi^k(b)$ converges weakly in $B$ for every $b\in B$.
\item[(c)] $\varphi$ is weakly mean ergodic.
\item[(d)] There is a UCP idempotent $\psi$ of $B$ onto $B^\varphi$ such that $\pi\circ\psi=\pi$ for every $\varphi$-harmonic representation $\pi$ of $B$.
\end{itemize}
In this case, $\psi=\rho_B$.
\end{proposition}
\begin{proof}
(a)$\Rightarrow$(b): Since $B\subseteq C^*(M^\varphi)$, Proposition~\ref{prop:sot} (Prunaru) gives $\varphi^k(b)\to\rho(b)=\rho_B(b)\in B$ in the weak$^*$ topology of $B^{**}$, and hence in the weak topology of $B$.
 
(b)$\Rightarrow$(c) is trivial. (c)$\Rightarrow$(d) is Lemma~\ref{lem:meanergodic}, and (c)$\Rightarrow$(a) follows from $\rho_B(b)=\psi_\omega(b)=\lim_{m\to\omega}\mathcal{C}_m(b)=\lim_{m\to\infty}\mathcal{C}_m(b)=\psi(b)\in B$.
 
(d)$\Rightarrow$(a): Let $W$ be the $\sigma$-weak closure of $\rho_B(B)$ in $H^\infty(M,\varphi)$ and $\sigma$ a faithful normal representation of the von Neumann algebra $W$. Then $\pi:=\sigma\circ\rho_B$ is a representation of $B$, and it is $\varphi$-harmonic because $\rho_B\circ\varphi=\rho_B$ (as $\psi_\omega\circ\varphi=\psi_\omega$). By (d), $\sigma(\rho_B(\psi(b)))=\sigma(\rho_B(b))$ for every $b\in B$. Since $\psi(b)\in B^\varphi\subseteq M^\varphi$ is fixed by $\psi_\omega$, we have $\rho_B(\psi(b))=\psi(b)$, so that both $\psi(b)$ and $\rho_B(b)$ lie in $W$, and the faithfulness of $\sigma$ on $W$ gives $\rho_B(b)=\psi(b)\in B$.
\end{proof} 
 
\begin{remark}[Concrete versus universal: Toeplitz operators with $L^\infty$-symbols]
\label{rem:Linfty}
Proposition~\ref{prop:rho_B} concerns the universal $W^*$-envelope; in a concrete representation the condition $\rho_B(B)\subseteq B$ is weaker than the weak mean ergodicity. Let $M=B(H^2)$, $\varphi=\mathrm{Ad}\,s^*$, and $B=C^*(M^\varphi)=C^*(\{T_f:f\in L^\infty(\mathbb T)\})$. Then $\rho(B)=(M^\varphi,\circ)=M^\varphi \subseteq B$ holds automatically (Remark~\ref{rem:toeplitz}). Nevertheless, $\varphi$ is not weakly mean ergodic on $B$. 
Indeed, let $u(z)=\exp\bigl(-\tfrac{1+z}{1-z}\bigr)$ and $x=T_uT_u^*$, the projection onto $uH^2$. Then $x\in B$ and $\rho(x)=T_u\circ T_{\bar u}=T_{|u|^2}=1$, so that $\varphi^k(x)\to1$ strongly by Proposition~\ref{prop:sot}. 
Let $\hat k_r(z)=\sqrt{1-r^2}\,(1-rz)^{-1}$ be the normalized reproducing kernel of $H^2$ at $r\in(0,1)$. Since $s^*\hat k_r=r\hat k_r$ and $s$ is an isometry, $\Vert s\hat k_r-\hat k_r\Vert^2=2-2r\to0$ as $r\to1$, so that every weak$^*$ cluster point $\chi$ of the states $(\,\cdot\,\hat k_r|\hat k_r)$ is $\varphi$-invariant, whereas $(x\hat k_r|\hat k_r)=|u(r)|^2\to0$; thus $\chi(x)=0\neq1=\chi(1)$. 
If $\varphi^k(x)$ converged weakly in $B$, its limit would be $1$ and $\chi(\varphi^k(x))=\chi(x)=0$ would converge to $\chi(1)=1$, a contradiction. Hence, condition (b) of Proposition~\ref{prop:rho_B} fails for $(B,\varphi)$: in the universal $W^*$-envelope $B^{**}$ one has $\rho_B(B)\not\subseteq B$, although $\rho(B)\subseteq B$ in the concrete one. 
Let $\mathcal B_0$ be the set of $b\in B$ whose Ces\`aro means converge in norm; by the proof of Lemma~\ref{lem:meanergodic}, it is a closed subspace equal to $B^\varphi\oplus\overline{\mathrm{ran}(\mathrm{id}-\varphi)}$. It contains $B^\varphi$ and the compact operators (as $s$ is pure), but not $x$. Since $\varphi(c)-c=s^*[c,s]$ and $[s,T_f]$ has rank at most one, $\varphi$ induces the identity on $B/\mathbb K$. It would be interesting to describe $\mathcal B_0$ in terms of the results of Davidson \cite{Dav77} on operators commuting with Toeplitz operators modulo the compact operators.
\end{remark}
 
We then prove the converse of Theorem~\ref{thm:bdry} under the equivalent conditions in Proposition~\ref{prop:rho_B}, obtain a Coburn-type exact sequence, describe where the non-commutative Poisson boundary sits relative to the C$^*$-envelope, and work out the examples (Toeplitz, sphere, Cuntz--Toeplitz and the non-commutative disc algebra, Toeplitz--Pimsner algebras, free groups).
 
\begin{theorem}[The converse for $S=B^\varphi$]
\label{thm:conv}
Let $B$ be a unital $C^*$-algebra, $\varphi$ be a UCP selfmap on $B$, $S=B^\varphi$ with $C^*(S)=B$. Suppose there is a UCP idempotent $\psi:B\to S$ such that $\pi\circ\psi=\pi$ for every representation $\pi$ with $\pi\circ\varphi=\pi$.
Then an irreducible representation $\pi$ of $B$ is a boundary representation for $S$ if and only if $\pi\circ\varphi=\pi$.
\end{theorem}
\begin{proof}
($\Rightarrow$) is Theorem~\ref{thm:bdry}(1).
 
($\Leftarrow$) Let $\theta:B\to B(H_\pi)$ be a UCP extension of $\pi\!\upharpoonright_S$.
Since $\psi(b)\in S$, $\theta(\psi(b))=\pi(\psi(b))=\pi(b)$ by the assumption, and hence $\theta\circ\psi=\pi$.
For each $x\in S$ put $d=\psi(x^*x)-x^*x\geq0$, which holds by the Schwarz inequality $\psi(x^*x)\ge\psi(x)^*\psi(x)=x^*x$.
Then $\theta(x^*x)=\theta(\psi(x^*x))-\theta(d)\leq\pi(x^*x)$, while the Schwarz inequality gives $\theta(x^*x)\geq\theta(x)^*\theta(x)=\pi(x)^*\pi(x)=\pi(x^*x)$.
Hence $\theta(x^*x)=\theta(x)^*\theta(x)$. Applying this to $x^*\in S$ as well, $x$ falls in $\mathrm{mult}(\theta)$. Thus $S\subseteq\mathrm{mult}(\theta)$, a C$^*$-algebra, whence $B=C^*(S)\subseteq\mathrm{mult}(\theta)$ and $\theta$ is a $*$-homomorphism agreeing with $\pi$ on $S$. Consequently, $\theta=\pi$ on $C^*(S)=B$.
\end{proof}
 
\begin{remark}
By Proposition~\ref{prop:rho_B}, the assumption of Theorem~\ref{thm:conv} is equivalent to the weak mean ergodicity of $\varphi$ (when $M=B^{**}$). Consequently, a counterexample to the converse with $S=B^\varphi$ must admit no UCP idempotent $\psi$ of $B$ onto $B^\varphi$ such that $\pi\circ\psi=\pi$ for every representation $\pi$ with $\pi\circ\varphi=\pi$; in particular, $\varphi$ is not weakly mean ergodic on $B$ (in the commutative case, there is a nonzero invariant signed measure annihilating $B^\varphi$). 
We do not know whether the implication $\varphi$-harmonicity $\Rightarrow$ being a boundary representation holds for $S=B^\varphi$ without mean ergodicity; any proof would have to avoid the UCP idempotent $\psi$ of Lemma~\ref{lem:meanergodic}. A natural test case is the algebra $B$ of Remark~\ref{rem:Linfty}: by Theorem~\ref{thm:choquet}, the boundary representations for $B^\varphi=\{T_f:f\in L^\infty(\mathbb T)\}$ are the representations $\chi\circ\rho$ with $\chi$ a character of $L^\infty(\mathbb T)$, and it is not clear whether every $\varphi$-harmonic irreducible representation of $B$ is of this form.
\end{remark}
 
The proof of the next corollary follows the same pattern as the construction of non-commutative Poisson boundaries in Section~\ref{sec:env}.
 
\begin{corollary}[Coburn-type exact sequence; the boundary is the spectrum]
\label{cor:coburn}
Assume $\varphi$ is weakly mean ergodic on $B$, $S=B^\varphi$, $C^*(S)=B$, and let $\psi$ be the limit of the Ces\`aro means (see Lemma~\ref{lem:meanergodic}). Set $x\circ y:=\psi(xy)$ on $S$. Then:
\begin{itemize}
\item[(1)] $(S,\circ)$ is a C$^*$-algebra, $\psi:B\to(S,\circ)$ is a unital surjective $*$-homomorphism with $\psi\!\upharpoonright_S=\mathrm{id}$, and $\ker\psi=J_S$. Hence $C^*_{\mathrm{env}}(S)\cong(S,\circ)$ and
\[
0\to J_S\to B\xrightarrow{\;\psi\;}C^*_{\mathrm{env}}(S)\to0
\]
is exact, with the unital completely positive section $S\hookrightarrow B$.
\item[(2)] The boundary representations for $S$ are exactly the representations $\hat\pi\circ\psi$, $\hat\pi$ an irreducible representation of $C^*_{\mathrm{env}}(S)$. In particular, every irreducible representation of the $C^*$-envelope is a boundary representation for $S$.
\end{itemize}
\end{corollary}
\begin{proof}
(1) By the Choi--Effros theorem \cite[Theorem 3.1]{CE1977} (see also \cite[Chapter~15]{Paulsen2002}), $(S,\circ)$ is a unital C$^*$-algebra and the identity map $\mathrm{id}:S\to(S,\circ)$ is a unital complete order isomorphism (hence a complete isometry).
Thus $\psi:B\to(S,\circ)$ is UCP.
For every $x\in S$, we have $\psi(x^*x)=x^*\circ x=\psi(x)^*\circ\psi(x)$, and likewise for $xx^*$.
Thus $S\subseteq\mathrm{mult}(\psi)$, a C$^*$-subalgebra, whence $B=C^*(S)\subseteq\mathrm{mult}(\psi)$ and $\psi$ is a $*$-homomorphism as a UCP map into $(S,\circ)$.
Since $\psi\!\upharpoonright_S=\mathrm{id}$ is completely isometric, $\ker\psi$ is a boundary ideal, and hence $\ker\psi\subseteq J_S$.
The quotient map $q:B\to B/J_S$ factors as $q=\sigma\circ\psi$ with $\sigma$ a $*$-homomorphism of $(S,\circ)$, isometric on $S$ (as so is $q\!\upharpoonright_S$), hence injective. Therefore $J_S=\ker q=\ker\psi$ and $B/J_S\cong(S,\circ)$.
 
(2) A harmonic $\pi$ satisfies $\pi\circ\psi=\pi$, so that it vanishes on $\ker\psi=J_S$ and is $\hat\pi\circ\psi$ with $\hat\pi$ irreducible. Conversely, $\psi\circ\varphi=\psi$ (Lemma~\ref{lem:meanergodic}), and hence $b-\varphi(b)\in\ker\psi$, and $\hat\pi\circ\psi$ is harmonic. Now the desired assertion follows from Theorem~\ref{thm:conv}.
\end{proof}
 
\begin{proposition}[Where the non-commutative Poisson boundary sits: atoms and the diffuse part]
\label{prop:poisson}
In the situation of Corollary~\ref{cor:coburn} put $A:=C^*_{\mathrm{env}}(S)=(S,\circ)$ and $(M,\varphi):=(B^{**},\varphi^{**})$. Then:
\begin{itemize}
\item[(1)] $\psi_\omega\!\upharpoonright_B=\psi$ for every $\omega$; hence $(S,\circ)$ is a C$^*$-subalgebra of the non-commutative Poisson boundary $H^\infty(M,\varphi)$ (the inclusion $\iota:A\hookrightarrow H^\infty(M,\varphi)$ is injective) and $\rho_B=\iota\circ\psi$.
\item[(2)] \textbf{Every} representation $\pi$ of $A$ (irreducible or not) is of the form $\pi=\hat\pi\circ\iota$ with $\hat\pi$ a normal representation of $H^\infty(M,\varphi)$. Consequently, the normal extension $\bar\iota:A^{**}\to H^\infty(M,\varphi)$ of $\iota$ is \emph{injective}, and the $\sigma$-weak closure $W$ of $\iota(A)$, which equals $\bar\rho_B(M)$, is isomorphic to $A^{**}$.
\item[(3)] Hence,
\[
W\;\cong\;\bigoplus_{[\pi]\in\hat A}B(\mathcal{H}_\pi)\ \oplus\ W_{\mathrm{na}},\qquad W_{\mathrm{na}}=(A^{**})_{\mathrm{na}}\ (\text{the non-atomic part}),
\]
the type I summands being exactly the summands $z_\pi\circ H^\infty(M,\varphi)\cong B(\mathcal{H}_\pi)$ of Theorem~\ref{thm:bdry} (one for each boundary representation of $S$), and for every representation $\pi$ of $A$, $\pi(A)''$ is isomorphic to a direct summand of $W$.
\end{itemize}
\end{proposition}
\begin{proof}
(1) Since $\varphi^k(b)\to\psi(b)$ weakly in $B$ (Lemma~\ref{lem:meanergodic}), i.e., weak$^*$ in $B^{**}$, we have $\psi_\omega(b)=\psi(b)$. For any $x,y\in S$, $xy\in B$ and $x\circ y=\psi_\omega(xy)=\psi(xy)$.
 
(2) Let $\pi$ be a representation of $A$. Since $\psi\circ\varphi=\psi$, $\pi\circ\psi$ is a representation of $B$ with $\pi\circ\psi\circ\varphi=\pi\circ\psi$, so that by Proposition~\ref{prop:harm}, $\pi\circ\psi=\hat\pi\circ\rho_B=\hat\pi\circ\iota\circ\psi$ with $\hat\pi$ normal, and $\psi$ is surjective, so $\pi=\hat\pi\circ\iota$. Then $\bar\pi=\hat\pi\circ\bar\iota$ on $A^{**}$ (both normal and equal on $A$). If $\pi$ is the universal representation, $\bar\pi$ is injective, hence so is $\bar\iota$.
Moreover, $\bar\iota(A^{**})$ is $\sigma$-weakly closed, contains $\iota(A)$ and is contained in its $\sigma$-weak closure; hence $\bar\iota(A^{**})=W$. Since $\rho_B(B)=\iota(\psi(B))=\iota(A)$, $W$ is also the $\sigma$-weak closure of $\rho_B(B)$, i.e., $W=\bar\rho_B(M)$.
 
(3) The atomic part of $A^{**}$ is $\bigoplus_{[\pi]}B(\mathcal{H}_\pi)$, and a normal image of $A^{**}$ is isomorphic to a direct summand of $A^{**}\cong W$. For an irreducible $\pi$, $\pi\circ\psi$ is a boundary representation for $S$ (Corollary~\ref{cor:coburn}), and the central support of $\hat\pi$ is the summand $z_\pi\circ H^\infty(M,\varphi)$ of Theorem~\ref{thm:bdry}.
Since $\hat\pi(W)=B(\mathcal H_\pi)$ (as $\hat\pi(\iota(A))=\pi(A)$ is $\sigma$-weakly dense and $\hat\pi$ is normal) and $\hat\pi$ maps $z_\pi\circ H^\infty(M,\varphi)$ isomorphically onto $B(\mathcal H_\pi)$, we have $z_\pi\circ W=z_\pi\circ H^\infty(M,\varphi)$; thus the type I summand of $W$ corresponding to $[\pi]$ is $z_\pi\circ H^\infty(M,\varphi)$. \end{proof}
 
\begin{remark}[Dictionary]
\label{rem:dictionary}
Boundary representations (maximal dilations, the non-commutative Choquet boundary) are the \emph{atoms} of the non-commutative Poisson boundary computed in the universal representation, that is, in $H^\infty(B^{**},\varphi^{**})$ (Proposition~\ref{prop:poisson}).
If $B\subseteq B(\mathcal{H})$ with $B''=N$ and $\varphi$ extends to a normal UCP map $\theta$ on $N$, then Proposition~\ref{prop:bridge}(2) (with $\sigma$ the inclusion) shows that $H^\infty(N,\theta)$ is a corner of $H^\infty(B^{**},\varphi^{**})$ by a central projection. In the examples below this quotient kills the atoms (the identity representation $\pi$ of $B$ does not annihilate $J_S$, so it is not harmonic; $H^\infty(N,\theta)$ has no normal irreducible representation): for the Toeplitz algebra $H^\infty(N,\theta)=L^\infty(\mathbb T)$ is the weak$^*$ closure of $A=C(\mathbb T)$, the characters $\mathrm{ev}_t$ (Dirac measures) are the boundary representations, and Lebesgue measure is the diffuse part; for the Cuntz--Toeplitz algebra $H^\infty(N,\theta)=\mathcal O_d''$ is the factor of \cite{BBDR}, the weak closure of $A=\mathcal O_d$, and boundary representations are the irreducible representations of $\mathcal O_d$.
Thus the non-commutative Poisson boundary of a concrete representation records the diffuse part, whereas the Choquet boundary records the atoms; this is a kind of Lebesgue-type decomposition (compare \cite{CD2016}).
\end{remark}
 
\begin{example}[Toeplitz operators]
Let $B=\mathcal T=C^*(s)=\{T_f+K:f\in C(\mathbb T),\ K\ \text{compact}\}\subset B(H^2)$ be the Toeplitz algebra and $\varphi(a)=s^*as$.
Then $\varphi^k(T_f+K)=T_f+s^{*k}Ks^k\to T_f$ in norm, because $s^{*k}Ks^k\to0$ in norm for every compact $K$ (the backward shift is pure).
Hence $\varphi$ is weakly mean ergodic, $B^\varphi=\{T_f:f\in C(\mathbb T)\}$ generates $\mathcal T$, and $\psi(T_f+K)=T_f$ has kernel $\mathbb K$, the ideal of compact operators on $H^2$.
By Theorem~\ref{thm:conv} (and Corollary~\ref{cor:coburn}), the boundary representations for $B^\varphi$ are exactly the harmonic irreducible representations, that is, the characters $\mathrm{ev}_t\circ q$, $t\in\mathbb T$, where $q:\mathcal T\to C(\mathbb T)$ is the symbol map (the identity representation is not harmonic). Corollary~\ref{cor:coburn} recovers Coburn's exact sequence
\[
0\to\mathbb K\to\mathcal T\xrightarrow{\;q\;}C(\mathbb T)\to0
\]
with $C^*_{\mathrm{env}}(B^\varphi)=C(\mathbb T)$.
\end{example}
 
\begin{example}[Toeplitz operators on the sphere]
Let $H^2$ be the Hardy space of the sphere $\mathbb S=S^{2d-1}\subset\mathbb{C}^d$, $T_f$ the Toeplitz operators, $B=C^*(\{T_f:f\in C(\mathbb S)\})=\{T_f+K:f\in C(\mathbb S),\ K\ \text{compact}\}$ (Coburn \cite{Coburn73}, Venugopalkrishna \cite{Venugopalkrishna72}), and $s_i=T_{z_i}$ with coordinate function $z_i$. Since $s_i^*s_i=T_{|z_i|^2}$ and $\sum_i|z_i|^2=1$,
\[
\varphi(b)=\sum_{i=1}^d s_i^*bs_i,\qquad b\in B,
\]
is UCP. For every $f\in C(\mathbb S)$, we have
\[
s_i^*T_fs_i=T_{\bar z_i}T_{fz_i}=T_{|z_i|^2f},
\]
implying $\varphi(T_f)=T_f$.
If $K=|\xi\rangle\langle\eta|$ is a rank-one operator with polynomials $\xi,\eta$ of degree less than $k$, then $\varphi^k(K)=0$, since each monomial in $s_i^*$ of length $k$ annihilates such polynomials (the monomials are orthogonal). Such operators span a norm-dense subspace of the compact operators and $\Vert \varphi^k\Vert\leq1$, and hence $\varphi^k(K)\to0$ for all compact $K$.
Therefore $B^\varphi=\{T_f:f\in C(\mathbb S)\}$, $\varphi$ is weakly mean ergodic, $\psi(T_f+K)=T_f$ for all $f\in C(\mathbb S)$ and compact $K \in \mathbb K$, and Corollary~\ref{cor:coburn} gives Coburn's sequence
\[
0\to\mathbb K\to B\xrightarrow{\;q\;}C(\mathbb S)\to0
\]
with $C^*_{\mathrm{env}}(\{T_f:f\in C(\mathbb S)\})=C(\mathbb S)$ (here $T_f\circ T_g=T_{fg}$, since the semicommutators are compact); the boundary representations for $\{T_f:f\in C(\mathbb S)\}$ are the characters $\mathrm{ev}_t\circ q$, $t\in\mathbb S$.
\end{example}
 
\begin{remark}[Drury--Arveson $d$-shift]
\label{rem:da}
Let $\alpha\in\mathbb N^d$ be a multi-index with $|\alpha|=\sum_j\alpha_j$, $\varepsilon_i\in\mathbb N^d$ the $i$-th unit multi-index, and $(e_\alpha)$ the orthonormal basis of monomials of a space of holomorphic functions on the ball. The $d$-shift $T=(T_1,\dots,T_d)$ on the Drury--Arveson space is $T_ie_\alpha=\sqrt{(\alpha_i+1)/(|\alpha|+1)}\,e_{\alpha+\varepsilon_i}$ \cite{Arv1998}, so that $\Lambda:=\sum_iT_i^*T_i=(N+d)/(N+1)$ on the level $|\alpha|=N$ is not $1$, and the map $\sum_iT_i^*\cdot T_i$ is not unital.
Normalizing, $R_i:=T_i\Lambda^{-1/2}$ satisfies $\sum_iR_i^*R_i=1$ and $R_ie_\alpha=\sqrt{(\alpha_i+1)/(|\alpha|+d)}\,e_{\alpha+\varepsilon_i}$, which are exactly the weights of the shifts $s_i=T_{z_i}$ on the Hardy space of the sphere. Hence $\varphi(a)=\sum_iR_i^*aR_i$ is unitarily equivalent to the map in the preceding example. Since $\Lambda-1$ is compact and $\Lambda^{\pm1/2}$ is a continuous function of $\sum_iR_iR_i^*=N/(N+d-1)$, we have $C^*(T)=C^*(R)\supseteq\mathbb K$, the compact operators, and Corollary~\ref{cor:coburn} recovers Arveson's extension 
\[
0\to\mathbb K\to C^*(T)\to C(\partial\mathbb B_d)\to0, 
\] 
see \cite{Arv1998}. Note that $B^\varphi$ is the image of the Hardy Toeplitz operators, not the operator system $\mathrm{span}\{1,T_i,T_i^*\}$. 
\end{remark}
 
\begin{example}[Cuntz--Toeplitz algebra and the non-commutative disc algebra]
\label{ex:ct}
Let $\mathcal F$ be the full Fock space over $\mathbb C^d$, whose standard basis consists of the vacuum vector $\Omega$ and all $e_\alpha=e_{\alpha_1}\otimes\cdots\otimes e_{\alpha_n}$ with multi-index $\alpha=(\alpha_1,\dots,\alpha_n)$. Let $l_i$ and $r_i$ be the left and right creation operators (they commute), $B=\mathcal E_d=C^*(l_1,\dots,l_d)=\overline{\mathrm{span}}\{l_\mu l_\nu^*\}\supseteq\mathbb K:=\mathbb K(\mathcal F)$, the compact operators on $\mathcal F$, where $l_\mu=l_{\mu_1}\cdots l_{\mu_n}$ with multi-index $\mu=(\mu_1,\dots,\mu_n)$. Then $\mathcal E_d/\mathbb K=\mathcal O_d$ (simple). Let $q:\mathcal E_d\to\mathcal O_d$ denote the quotient map, and
\[
\varphi(b)=\tfrac1d\sum_{i=1}^dr_i^*br_i,\qquad b\in B
\]
is a UCP map (Kraus operators $r_i/\sqrt d$). Write $\theta_{\alpha,\beta}=|e_\alpha\rangle\langle e_\beta|$ and let $\mu^-$ denote $\mu$ without its last letter $\mu_{\mathrm{last}}$. A direct computation (using $r_i^*l_\mu r_j=\delta_{ij}l_\mu$ and $l_\nu^*r_j=r_jl_\nu^*+\delta_{j,\nu_{\mathrm{last}}}|\Omega\rangle\langle e_{\nu^-}|$) gives
\[
\varphi(l_\mu l_\nu^*)=l_\mu l_\nu^*+\tfrac1d\,\delta_{\mu_{\mathrm{last}},\nu_{\mathrm{last}}}\,\theta_{\mu^-,\nu^-}\qquad(|\mu|,|\nu|\geq1),
\]
and $\varphi(l_\mu)=l_\mu$, $\varphi(l_\nu^*)=l_\nu^*$; moreover $\varphi(\theta_{\alpha,\beta})=\frac1d\delta_{\alpha_{\mathrm{last}},\beta_{\mathrm{last}}}\theta_{\alpha^-,\beta^-}$ (and $0$ if $\alpha$ or $\beta$ is empty).
Hence $\varphi(\mathcal E_d)\subseteq\mathcal E_d$, $\varphi(b)-b\in\mathbb K$ for all $b\in\mathcal E_d$ (so that $\varphi$ induces the identity on $\mathcal O_d$), and $\varphi^n(l_\mu l_\nu^*)$ is eventually constant for every word (the correction terms are matrix units, which are annihilated after finitely many steps; $\mu^{(r)}$ denotes $\mu$ without its last $r$ letters). Since the words $l_\mu l_\nu^*$ span a dense subspace of $\mathcal E_d$ (this includes the matrix units, so that $\varphi^n\to0$ on $\mathbb K$) and $\Vert\varphi^n\Vert\leq1$, $\varphi^n(b)$ converges in norm for every $b\in\mathcal E_d$, and
\[
\psi(l_\mu l_\nu^*)=l_\mu l_\nu^*+\sum_{r\ge1}d^{-r}\,[\text{last $r$ letters of $\mu,\nu$ agree}]\;\theta_{\mu^{(r)},\nu^{(r)}}.
\]
Hence $\varphi$ is weakly mean ergodic, $S:=B^\varphi=\psi(\mathcal E_d)$ contains $\mathcal A_d+\mathcal A_d^*$ ($\mathcal A_d$ the non-commutative disc algebra), $S\cap\mathbb K=0$, and $C^*(S)=\mathcal E_d$.
Corollary~\ref{cor:coburn} gives $(S,\circ)\cong\mathcal O_d$ (via $q\!\upharpoonright_S$), $J_S=\mathbb K$, the exact sequence
\[
0\to\mathbb K\to\mathcal E_d\overset{\,q\,}{\to}\mathcal O_d\to0,
\]
and the boundary representations for $S$ are exactly $\hat\pi\circ\psi$ with irreducible representation $\hat\pi$ of $\mathcal O_d$. This is consistent with $C^*_{\mathrm{env}}(\mathcal A_d)=\mathcal O_d$.
 
A mirror-image situation on $B(\mathcal F)$ was studied by Bhat, Bikram, De and Rakshit \cite{BBDR}: for the Markov operator built from the left creation operators, they showed that the non-commutative Poisson boundary is the $\sigma$-weak closure of the Cuntz algebra generated by the right creation operators $r_i$, and that the Choi--Effros product satisfies $r_i\circ r_j^*=r_ir_j^*+\frac1d\delta_{ij}|\Omega\rangle\langle\Omega|$ \cite[Proposition~5]{BBDR} (the unitary reversing words exchanges $l_i$ and $r_i$, and our $\psi(l_\mu l_\nu^*)$ corresponds to $r_\mu\circ r_\nu^*$).
\end{example}

\begin{proposition}[Boundary representations of the free disc system]
\label{prop:disc}
In Example~\ref{ex:ct}, let $S_d=\overline{\mathrm{span}}\{1,l_\mu,l_\nu^*\}=\mathcal A_d+\mathcal A_d^*\subset\mathcal E_d$ be the free disc operator system ($\mathcal A_d$ the non-commutative disc algebra). Then $S_d\subseteq T:=B^\varphi$ and $\partial S_d=\partial T=\mathrm{Irr}(\mathcal O_d)$. More precisely, every representation $\pi$ of $\mathcal O_d$ (irreducible or not) is the unique UCP extension to $\mathcal E_d$ of $\pi\!\upharpoonright_{S_d}$ (here $\pi$ is regarded as a representation of $\mathcal E_d$ through $q$).
\end{proposition}
\begin{proof}
Let $\pi$ be a representation of $\mathcal O_d$, $L_i=\pi(l_i)$, so $L_i^*L_i=1$ and $\sum_iL_iL_i^*=1$, and let $\theta:\mathcal E_d\to B(\mathcal H)$ be a UCP map with $\theta\!\upharpoonright_{S_d}=\pi\!\upharpoonright_{S_d}$. By the Schwarz inequality, $\theta(l_il_i^*)\ge L_iL_i^*$, while $\sum_i\theta(l_il_i^*)=1-\theta(|\Omega\rangle\langle\Omega|)\leq1=\sum_iL_iL_i^*$; hence $\theta(l_il_i^*)=L_iL_i^*$. Together with $\theta(l_i^*l_i)=1=L_i^*L_i$, this gives $l_i\in\mathrm{mult}(\theta)$, so that $\mathcal E_d=C^*(l_i)\subseteq\mathrm{mult}(\theta)$ and $\theta$ is a $*$-homomorphism with $\theta(l_i)=L_i$, i.e., $\theta=\pi$. Hence $\mathrm{Irr}(\mathcal O_d)\subseteq\partial S_d$. Conversely, $\partial S_d\subseteq\partial T$ (extensions of $\pi\!\upharpoonright_T$ are extensions of $\pi\!\upharpoonright_{S_d}$) and $\partial T=\mathrm{Irr}(\mathcal O_d)$ by Corollary~\ref{cor:coburn}.
\end{proof}

\begin{remark}[The free Lebesgue measure: atoms and the diffuse part]
\label{rem:lebesgue}
For $d=1$, $S_1=\mathcal A(\mathbb D)+\overline{\mathcal A(\mathbb D)}$, the boundary representations are the Dirac masses on $\mathbb T$ (atoms of $C(\mathbb T)^{**}$), the non-commutative Poisson boundary of $(B(H^2),\mathrm{Ad}\,s^*)$ is $L^\infty(\mathbb T)$, and the F.\ and M.\ Riesz theorem says that the measures annihilating $\mathcal A(\mathbb D)$ are absolutely continuous, so that Lebesgue measure has a unique extension to $C(\mathbb T)$.
For $d\ge2$, the free Lebesgue measure $\chi(x)=\langle\lambda(x)\xi_\emptyset,\xi_\emptyset\rangle$ on $S_d$ (the left regular representation $\lambda$ of the free semigroup algebra; it is absolutely continuous in the sense of \cite{DLP05}) has non-equivalent extensions of both kinds: the vacuum state of the Poisson boundary $\mathcal O_d''$ of Example~\ref{ex:ct} restricts to $\chi$ (it is the state with $\omega(l_\mu l_\nu^*)=\delta_{\mu\nu}d^{-|\mu|}$, of type III$_{1/d}$ \cite{BBDR}), and so does the vector state of a unit wandering vector (i.e., $\{L_w\xi\}_w$ orthonormal) of any representation of $\mathcal O_d$; irreducible representations of $\mathcal O_d$ with a wandering vector exist (the atomic representations of \cite[Section~6]{DLP05}), and they are boundary representations of $\mathcal A_d$ by Proposition~\ref{prop:disc}. Thus the dichotomy atom versus diffuse part, which is canonical on $\mathcal O_d^{**}$, does not coincide with the absolutely continuous versus singular decomposition of functionals on $\mathcal A_d$ \cite{DLP05,Kennedy,JuryMartin}. 
\end{remark}
 
\begin{remark}[Weak$^*$ density fails in the universal setting]
\label{rem:notdense}
In Example~\ref{ex:ct} ($B=\mathcal E_d$, $M=B^{**}$) the $\sigma$-weak closure $W$ of $S=B^\varphi$ (which is $\bar\rho_B(M)$ by Proposition~\ref{prop:poisson}(2)) is strictly smaller than $H^\infty(M,\varphi)$, and $\psi_\omega$ is not normal, although $\varphi$ is norm mean ergodic on $B$.
Indeed, let $p_0\in M$ be the central support of $\mathbb K$, so that $\bar\pi_0\!\upharpoonright_{p_0M}:p_0M\to B(\mathcal F)$ is an isomorphism, where $\bar\pi_0$ is the normal extension of the identity representation, and $\bar\pi_0\circ\varphi=\theta\circ\bar\pi_0$ with the normal UCP map $\theta=\frac1d\sum_ir_i^*\,\cdot\,r_i$ on $B(\mathcal F)$.
Since $\varphi(\mathbb K)\subseteq\mathbb K$, we have $\varphi(p_0)\in p_0M$, and hence $\varphi(p_0)\leq p_0$, and $\bar\pi_0(\varphi(p_0))=\theta(1)=1=\bar\pi_0(p_0)$ gives $\varphi(p_0)=p_0$; thus $p_0\in M^\varphi$.
 
(i) If $(e_k)$ is an approximate unit of $\mathbb K$, then $\psi_\omega(e_k)=\psi(e_k)=0$ and $e_k\nearrow p_0$; normality of $\psi_\omega$ would give $\psi_\omega(p_0)=0\neq p_0$.
 
(ii) The vector state $g=(\bar\pi_0(\cdot)\Omega|\Omega)$ is normal with $g(p_0)=1$, and $g\!\upharpoonright_S=h\circ q\!\upharpoonright_S$ for a state $h$ of $\mathcal O_d$, because $q\!\upharpoonright_S:S\to\mathcal O_d$ is a unital complete order isomorphism (Corollary~\ref{cor:coburn}). If $p_0\in W$, there is a net $t_\alpha\in S$ with $t_\alpha\to p_0$ weak$^*$; then $g(t_\alpha)\to1$, whereas $q(t_\alpha)=\bar q(t_\alpha)\to\bar q(p_0)=0$ $\sigma$-weakly (as $\bar q$ kills $p_0$), so that $g(t_\alpha)=h(q(t_\alpha))\to0$, a contradiction. Hence $p_0\in M^\varphi\setminus W$.
\end{remark}

\begin{example}[Toeplitz--Pimsner algebra over $C(X)$]
\label{ex:pim}
Let $X$ be a compact space, $\tau$ a homeomorphism of $X$ and $\mu$ a $\tau$-invariant probability measure with full support, and let $E=\mathbb C^d\otimes C(X)$ be the $C(X)$-correspondence whose left action is twisted by $\tau^{-1}$ as in \cite[Example~4.3]{Izumi2016}. Let $\mathcal W$ be the set of words in $\{1,\dots,d\}$, $P_k$ the projection of $\ell^2(\mathcal W)$ onto the words of length $k$, and $l_i$, $r_i$ the left and right creation operators. On $H_\mu=\ell^2(\mathcal W)\otimes L^2(X,\mu)$ we define
\[
\Pi(a)=\sum_kP_k\otimes M_{a\circ\tau^{-k}},\qquad T_i=l_i\otimes1,\qquad\rho_i=r_i\otimes U,\qquad Uf=f\circ\tau^{-1}
\]
($U$ is unitary since $\mu$ is invariant), let $B=C^*(\Pi(C(X)),T_1,\dots,T_d)$, which is the Toeplitz--Pimsner algebra $\mathcal T_E$ (the Fock representation induced from a faithful representation of $C(X)$ is faithful \cite{Pimsner}), let $\mathbb K_E$ be the ideal generated by $1-\sum_iT_iT_i^*$, so that $B/\mathbb K_E=\mathcal O_E$, and let $\varphi(x)=\frac1d\sum_i\rho_i^*x\rho_i$.
Since $\rho_i$ commutes with $\Pi(a)$ and $T_j$ (use $UM_gU^*=M_{g\circ\tau^{-1}}$ and that $l_j$, $r_i$ commute) and $\rho_i^*\rho_j=\delta_{ij}$, the map $\varphi$ is UCP and fixes $\Pi(a)$, $T_j$, $T_j^*$. The computation of Example~\ref{ex:ct} (using $l_v^*r_i=r_il_v^*+\delta_{i,v_{\mathrm{last}}}|\Omega\rangle\langle e_{v^-}|$) gives, for words $u,v$ of length $\ge1$,
\[
\varphi(T_u\Pi(a)T_v^*)=T_u\Pi(a)T_v^*+\tfrac1d\,\delta_{u_{\mathrm{last}},v_{\mathrm{last}}}\,\theta_{u^-,v^-}\otimes M_{a\circ\tau},
\]
and similarly for $\theta_{\alpha,\beta}\otimes M_g=T_\alpha\Pi(g)(1-\sum_iT_iT_i^*)T_\beta^*\in B$. Hence $\varphi(B)\subseteq B$, and on the dense span of the elements $T_u\Pi(a)T_v^*$ the iterates $\varphi^n$ are eventually constant (each application removes one letter), so that $\varphi^n(x)$ converges in norm for every $x\in B$ to some $\psi(x)$ with $x-\psi(x)\in\mathbb K_E$ and $\psi(\mathbb K_E)=0$. Thus $\ker \psi=\mathbb K_E$, $B^\varphi=\psi(B)$ contains the generators of $B$, and Corollary~\ref{cor:coburn} gives $C^*_{\mathrm{env}}(B^\varphi)=B/\mathbb K_E=\mathcal O_E$, whose irreducible representations are exactly the boundary representations for $B^\varphi$.
Only the invariance of $\mu$ is used, neither the continuity of $\tau$-orbits nor the structure of the support. If $X$ is a singleton, then we recover Example~\ref{ex:ct}. 
\end{example}

\begin{example}[Boundary action of free group]
\label{ex:free}
This example concerns a C$^*$-algebra rather than a specific operator system, but it exhibits the same picture as Coburn's sequence: the compact operators form the ideal, the quotient is the boundary crossed product, and the non-commutative Poisson boundary of the concrete representation is the weak closure of the quotient.
 
Let $\Gamma=\mathbb F_n$ with $n\geq2$, $\mu$ be a uniform probability measure on the $d=2n$ generators and inverses, $\bar\Gamma=\Gamma\cup\partial\Gamma$ the end-compactification ($\partial\Gamma$ is the ideal boundary, see e.g.\ \cite[Section 5.1]{BO08}) and
$B=C(\bar\Gamma)\rtimes_r\Gamma\subset B(\ell^2\Gamma)\supseteq C_0(\Gamma)\rtimes_r\Gamma=\mathbb K$ (the $C^*$-subalgebras of $B(\ell^2(\Gamma))$ generated by $C(\bar\Gamma)$ and $C_0(\Gamma)$, respectively, by multiplication and the left regular representation $\lambda$), so that
\[
0\to\mathbb K\to B\to C(\partial\Gamma)\rtimes_r\Gamma\to0
\]
is exact since $\Gamma$ is an exact group (see e.g.\ \cite[Section 5.1]{BO08}), and $C(\partial\Gamma)\rtimes_r\Gamma\cong\mathcal O_A$, the Cuntz--Krieger algebra with the $d\times d$ matrix $A$ determined by $A_{xy}=1$ unless $y=x^{-1}$; see Spielberg \cite{Spi91}.
Let $\varphi(x)=\sum_s\mu(s)\rho(s)x\rho(s)^*$, $x\in B(\ell^2(\Gamma))$, with the right regular representation $\rho$ (this is the Markov operator $Q_\mu$ of \cite{Izumi2004}; we also write $\varphi$ for its restriction to $B$).
Since $\rho$ commutes with $\lambda$, $\varphi(f\lambda(g))=(P_\mu f)\lambda(g)$ with the random walk operator $P_\mu$ defined by $P_\mu f(x)=\sum_s\mu(s)f(xs)$. Hence $\varphi(B)\subseteq B$ and $\varphi=\mathrm{id}$ modulo the compact operators $\mathbb K$ (right translations fix the points of $\partial\Gamma$, so $P_\mu f-f\in C_0(\Gamma)$ for $f\in C(\bar\Gamma)$).
For every $f\in C(\bar\Gamma)$, $P_\mu^nf\to Hf$ pointwise on $\bar\Gamma$ (we use standard facts on random walks on the free group: the walk converges almost surely to $\partial\Gamma$ and the Dirichlet problem at infinity is solvable \cite[Chapter~IV, Corollary 21.12]{Woess00}), where $Hf$ is the Poisson integral of $f\!\upharpoonright_{\partial\Gamma}$, thus weakly in $C(\bar\Gamma)$ (bounded pointwise convergence), and hence weakly in $B$ for $\varphi^n(f\lambda(g))=(P_\mu^nf)\lambda(g)$. By the proof of Lemma~\ref{lem:meanergodic}, the Ces\`aro means of $f\lambda(g)$ then converge in norm, and since the set of elements with norm convergent Ces\`aro means is a closed subspace, $\varphi$ is weakly mean ergodic with $\psi(f\lambda(g))=H(f\!\upharpoonright_{\partial\Gamma})\lambda(g)$.
 
Since Poisson extensions of members of $C(\partial\Gamma)$ separate the points of $\bar\Gamma$, they generate $C(\bar\Gamma)$ by the Stone--Weierstrass theorem; together with $\lambda(g)\in B^\varphi$ (as $\rho$ commutes with $\lambda$), we conclude that $C^*(B^\varphi)=B$.
(Here $\nu_x$ denotes the harmonic measure of the walk starting at $x\in\Gamma$ and $\nu=\nu_e$. Indeed, $Hf=f$ on $\partial\Gamma$; $\nu_x$ is not a Dirac mass; and $\nu_x\neq\nu_y$ for $x\neq y$ in $\Gamma$, because $\nu_g=g_*\nu_e$ and, if $a$ is the first letter of $g\neq e$ and $C_a$ the set of ends beginning with $a$, then $\nu_g(C_a)=(1-F)+F/d>1/d=\nu_e(C_a)$ (by the symmetry of $\mu$) with $F<1$ (transience) the probability of ever hitting $e$ from $g$.)
Moreover, $\mathbb K=C_0(\Gamma)\rtimes_r\Gamma\subseteq\ker\psi$ because $Hf=0$ for $f\in C_0(\Gamma)$, and $\ker\psi\neq B$ since $\psi(1)=1$; as $\mathcal O_A$ is simple \cite{Spi91}, $\ker\psi=\mathbb K$.
Now Corollary~\ref{cor:coburn} gives $C^*_{\mathrm{env}}(B^\varphi)\cong C(\partial\Gamma)\rtimes_r\Gamma=\mathcal O_A$, and the boundary representations are the irreducible representations of $\mathcal O_A$, and moreover $H^\infty(B(\ell^2\Gamma),Q_\mu)=L^\infty(\partial\Gamma,\nu)\rtimes\Gamma$, the weak closure of $\mathcal O_A$ in the harmonic-measure representation by \cite{Izumi2004}.
\end{example}

\appendix
\section{A $2$-positive variant of Arveson's unitary equivalence criterion}
\label{app:2pos}
 
This appendix is independent of the main body of the paper. We thank Yasuhiko Sato for pointing out in 2019 that Choi's theorem on multiplicative domains concerns $2$-positive maps, which led to this appendix.
 
\medskip
A unital linear bijection $\phi$ between operator systems is called a \emph{$k$-order isomorphism} if $\phi$ and $\phi^{-1}$ are $k$-positive. For a matrix $A\in M_n$ write $S_A=\mathrm{span}\{I,A,A^*\}$.
 
Arveson's unitary equivalence criterion \cite{Arv1970,Arv1972} says that if $A\in M_n$ is irreducible and $\phi\colon S_A\to S_B$ is a unital completely isometric bijection (equivalently, a complete order isomorphism) with $\phi(A)=B$, then $A$ and $B$ are unitarily equivalent. Its proof rests on two pillars: Arveson's extension theorem, used to extend $\phi$ and $\phi^{-1}$ to unital completely positive maps $\tilde\phi,\tilde\psi$ on $M_n$ (step~(1)), and the boundary theorem, which gives $\tilde\psi\circ\tilde\phi=\mathrm{id}$ (step~(2)); then the Schwarz inequality makes $\tilde\phi$ a $*$-isomorphism (step~(3)). The boundary theorem and the Schwarz argument work for $2$-positive maps (Proposition~\ref{prop:twoext}), whereas the extension theorem is not generalized beyond complete positivity: already unital positive maps, even unital positive isometries, on operator subsystems of matrix algebras need not have positive extensions \cite{CDPR}, and Proposition~\ref{thm:counter2}\textup{(iii)} below shows that the same happens for $2$-positive maps. The question is whether the unitary equivalence criterion survives without the extension theorem. We show that it does not for $2$-order isomorphisms (Proposition~\ref{thm:counter2}).
 
\subsection*{A.1. The boundary theorem for $2$-positive maps}
 
\begin{lemma}[Choi--Effros product for $2$-positive idempotents]
\label{lem:CE2}
Let $\varphi$ be a unital normal $2$-positive map on a von Neumann algebra $M$ and $\psi=\lim_{m\to\omega}\frac1{m+1}\sum_{k=0}^m\varphi^k$ ($\sigma$-weak limit along a free ultrafilter $\omega$), a unital $2$-positive idempotent onto $M^\varphi$ with $\varphi\circ\psi=\psi\circ\varphi=\psi$. Then $a\circ b:=\psi(ab)$ is an associative product on $M^\varphi$, $(M^\varphi,\circ)$ is a $W^*$-algebra, and the identity map $M^\varphi\to(M^\varphi,\circ)$ is a unital isometric order isomorphism.
\end{lemma}
\begin{proof}
For associativity it suffices to show $\psi(a^*x)=\psi(a^*\psi(x))$ and $\psi(xa)=\psi(\psi(x)a)$ for $a\in M^\varphi$ and $x=x^*\in M$. Apply the Kadison inequality for the unital positive map $\psi^{(2)}=\psi\otimes\mathrm{id}_{M_2}$ to $X=\left(\begin{smallmatrix}0&a^*\\a&x\end{smallmatrix}\right)$, and then $\psi^{(2)}$ to the resulting inequality (using $\psi\circ\psi=\psi$ and the positivity of $\psi^{(2)}$): we get
\[
\begin{pmatrix}\psi(a^*a)&\psi(a^*x)\\\psi(xa)&\psi(aa^*+x^2)\end{pmatrix}\ge\begin{pmatrix}\psi(a^*a)&\psi(a^*\psi(x))\\\psi(\psi(x)a)&\psi(aa^*+\psi(x)^2)\end{pmatrix},
\]
and the difference is a positive matrix with $(1,1)$-entry $0$, so its off-diagonal entries vanish. The Schwarz inequality for $\psi$ gives $\Vert a\Vert^2=\Vert\psi(a)^*\psi(a)\Vert\leq\Vert a^*\circ a\Vert\leq\Vert a^*a\Vert=\Vert a\Vert^2$, so that $(M^\varphi,\circ)$ is a C$^*$-algebra with the same norm, and it is a von Neumann algebra, since $M^\varphi$ is a dual Banach space. A unital isometry between operator systems is an order isomorphism, because positivity in a unital C$^*$-algebra is characterized by norms.
\end{proof}
 
\begin{proposition}
\label{prop:HU2pos}
Theorem~1 of \cite{HU2018} holds if the unital maps there are only assumed to be $2$-positive. They fail for unital positive maps: the transposition on $M_2(\mathbb C)$ is a counterexample.
\end{proposition}
\begin{proof}
The proofs in \cite{HU2018} use only the Schwarz inequality, Choi's multiplicative domains (both valid for $2$-positive maps \cite{Choi1974}) and the C$^*$-algebra structure of the non-commutative Poisson boundary, which is provided by Lemma~\ref{lem:CE2}. The transposition on $M_2(\mathbb C)$ fixes the irreducible set of symmetric matrices, but it is not the identity.
\end{proof}
 
\subsection*{A.2. Two-sided $2$-positive extensions give unitary equivalence}
 
\begin{proposition}
\label{prop:twoext}
For $i=1,2$ let $A_i\subseteq B(\mathcal H_i)$ be a unital subspace such that $C^*(A_i)$ is irreducible and the Calkin map $B(\mathcal H_i)\to B(\mathcal H_i)/K(\mathcal H_i)$ is not completely isometric on $A_i+A_i^*$ (this is automatic if $\dim\mathcal H_i<\infty$). Let $\pi\colon A_1\to A_2$ be a unital linear bijection such that $\pi$ and $\pi^{-1}$ admit unital $2$-positive extensions $\varphi\colon B(\mathcal H_1)\to B(\mathcal H_2)$ and $\psi\colon B(\mathcal H_2)\to B(\mathcal H_1)$. Then $\pi$ is implemented by a unitary $\mathcal H_1\to\mathcal H_2$.
\end{proposition}
\begin{proof}
First, $K(\mathcal H_i)\subseteq C^*(A_i)$: otherwise $C^*(A_i)\cap K(\mathcal H_i)=0$ by irreducibility, so the Calkin map would be an injective $*$-homomorphism on $C^*(A_i)$, hence completely isometric on $A_i+A_i^*$. The unital $2$-positive map $\psi\circ\varphi$ on $B(\mathcal H_1)$ fixes $A_1+A_1^*$, so that $C^*(B(\mathcal H_1)^{\psi\circ\varphi})\supseteq K(\mathcal H_1)$, a simple essential ideal, and the quotient map by $K(\mathcal H_1)$ is not completely isometric on $B(\mathcal H_1)^{\psi\circ\varphi}$. By Proposition~\ref{prop:HU2pos}, $\psi\circ\varphi=\mathrm{id}$, and by symmetry $\varphi\circ\psi=\mathrm{id}$.
For a unitary $u\in B(\mathcal H_1)$ the Schwarz inequality gives $1=\psi(\varphi(u))^* \psi(\varphi(u))\leq\psi(\varphi(u)^*\varphi(u))\leq\psi(\varphi(u^*u))=1$, so $\psi(d)=0$ for $d=1-\varphi(u)^*\varphi(u)\ge0$, and $d=\varphi\circ\psi(d)=0$; the same holds for $u^*$. Thus all unitaries lie in the multiplicative domain of $\varphi$, so that $\varphi$ is a $*$-homomorphism, bijective with inverse $\psi$, and hence implemented by a unitary \cite[Corollary~16.14]{Conway1999}.
\end{proof}
 
\subsection*{A.3. A $2$-order isomorphism need not be a unitary equivalence}

This subsection concerns operator systems in matrix algebras, so that, unlike in the rest of the paper, we denote matrices by capital letters.
 
\begin{lemma}
\label{lem:swap}
Let $a,b>0$, $E_0=E_0^*\in M_2(\mathbb C)$, $E_1\in M_2(\mathbb C)$ and
$T(a,b;E_0,E_1)=\begin{pmatrix}E_0&aE_1&0\\aE_1^*&E_0&bE_1\\0&bE_1^*&E_0\end{pmatrix}\in M_6(\mathbb C)$. Then $T(a,b;E_0,E_1)\ge0\iff T(b,a;E_0,E_1)\ge0$.
\end{lemma}
\begin{proof}
We may and do assume $E_0=I$ and $E_1=F$ (replace $E_0$ by $E_0+\varepsilon I$ and conjugate by $\mathrm{diag}(E_0^{-1/2},E_0^{-1/2},E_0^{-1/2})$). Then $T=I+K$, and conjugation by $\mathrm{diag}(1,-1,1)\otimes I$ changes $K$ into $-K$, so $T\ge0 \iff -I \le K \le I \iff \Vert K\Vert \leq1$. Conjugation by $\left[\begin{smallmatrix} 1 & 0 & 0 \\ 0 & 0 & 1 \\ 0 & 1 & 0 \end{smallmatrix}\right]$ changes $K$ into $\left(\begin{smallmatrix}0&W^* \\ W&0\end{smallmatrix}\right)$ with $W = [ a F^*, b F]$, and hence $\Vert K\Vert^2=\Vert a^2F^*F+b^2FF^*\Vert$. Write $P=F^*F$ and $Q=FF^*$; they are unitarily equivalent positive $2\times2$ matrices, so $P=tI+(s-t)p$ and $Q=tI+(s-t)q$ with rank-one projections $p,q$. If $p\ne q$, choose unit vectors $u\in p\mathbb C^2$, $w\in q\mathbb C^2$ with $\langle u,w\rangle\ge0$ and the unitary $V$ with $Vu=w$, $Vw=u$; then $VPV^*=Q$ and $VQV^*=P$, whence $\Vert a^2P+b^2Q\Vert=\Vert a^2Q+b^2P\Vert$, which is the same criterion for $T(b,a;I,F)$.
\end{proof}

It is well known, see e.g., \cite[Theorem 6.1]{Paulsen2002}, that any $n$-positive map into $M_n$ must be completely positive. We freely use this standard fact below. 
 
\begin{proposition}[A $2$-order isomorphism which is not a unitary equivalence]
\label{thm:counter2}
Let $a,b>0$, $a\ne b$, $A=J_3(a,b)=aE_{12}+bE_{23}$, $B=J_3(b,a)$, and let $\phi\colon S_A\to S_B$ be the linear map with $\phi(I)=I$, $\phi(A)=B$, $\phi(A^*)=B^*$ (so $\phi(h)=Ph^TP$ with $P=\mathrm{antidiag}(1,1,1)$). Then:
\begin{enumerate}
\item[\textup{(i)}] $A$ is irreducible and $A$ is not unitarily equivalent to $B$;
\item[\textup{(ii)}] $\phi$ and $\phi^{-1}$ are unital and $2$-positive, i.e., $\phi$ is a $2$-order isomorphism;
\item[\textup{(iii)}] at least one of $\phi$, $\phi^{-1}$ is not $3$-positive, and at least one of them has no unital $2$-positive extension to $M_3$; for $n=3$ the $3$-order version of Arveson's theorem holds.
\end{enumerate}
\end{proposition}
\begin{proof}
(i) $A$ is clearly irreducible, and hence $C^*(A)=M_3(\mathbb C)$. Traces of words in matrices are unitary invariants and $\mathrm{tr}(AAA^*AA^*A^*)=a^2b^4\ne a^4b^2=\mathrm{tr}(BBB^*BB^*B^*)$.
(ii) An element of $M_2(S_A)$ is self-adjoint if and only if it is $E_0\otimes I+E_1\otimes A+E_1^*\otimes A^*=T(a,b;E_0,E_1)$ for some $E_0=E_0^*, E_1 \in M_2(\mathbb C)$, and $\phi_2$ maps it to $T(b,a;E_0,E_1)$; positivity is preserved by Lemma~\ref{lem:swap}, and the same applies to $\phi^{-1}$.
(iii) If both were $3$-positive, they would be completely positive, and thus $\phi$ would be completely isometric and $A$ would be unitarily equivalent to $B$ by \cite{Arv1970,Arv1972}, contradicting (i). If both had unital $2$-positive extensions to $M_3$, Proposition~\ref{prop:twoext} would make $A$ unitarily equivalent to $B$.
\end{proof}
 
\begin{remark}[Where the obstruction is]
\label{rem:obstruction}
Thus, for a $2$-order isomorphism the extension step (1) of Arveson's proof may fail, whereas steps (2) and (3) survive (Proposition~\ref{prop:twoext}). Note that a unital $2$-order isomorphism is automatically an isometry, since $\Vert x\Vert\leq1$ if and only if $\left(\begin{smallmatrix}1&x\\x^*&1\end{smallmatrix}\right)\ge0$ in $M_2(S)$. 
Equality of the norms for \emph{scalar} coefficients does not imply unitary equivalence \cite[Proposition~3.1]{FGS11}. We do not know whether a unital $2$-isometric bijection $S_A\to S_B$ (that is, equality of $\Vert\alpha\otimes1+\beta\otimes A+\gamma\otimes A^*\Vert$ and $\Vert\alpha\otimes1+\beta\otimes B+\gamma\otimes B^*\Vert$ for all $\alpha,\beta,\gamma\in M_2(\mathbb C)$) forces $A$ to be unitarily equivalent to $B$. The natural approach through extensions would require a $2$-positive extension theorem for $2$-contractive maps, which we do not have, and a Hahn--Banach type argument cannot be expected to provide one: already unital contractive maps (hence positive ones) on operator subsystems of $M_{2n}(\mathbb C)$ need not have positive extensions, and the same holds for unital positive isometries of such systems \cite{CDPR}.
\end{remark}

\end{document}